\documentclass[aos,MSNbibl,dvips]{arximspdf}

\doi{10.1214/12-AOS1065} %kopijuoti is PTS
\volume{41}
\issue{1}
\pubyear{2013}
\firstpage{370}
\lastpage{400}

\makeatletter

\newproclaim{Remark}{Remark}

\newtheorem{theorem}{Theorem}
\newtheorem{proposition}{Proposition}
\newtheorem{lemma}{Lemma}

\newproclaim{definition}{Definition}
\newproclaim{example}{Example}

\newcommand{\E}{\mathbb{E}}
\newcommand{\Nat}{\mathbb{N}}
\newcommand{\indicator}{\mathbb{I}}

\newcommand{\conv}{\operatorname{conv}}
\newcommand{\real}{\mathbb{R}}
\newcommand{\Qcal}{\mathcal{Q}}

\newcommand{\ess}{\operatorname{ess}}
\newcommand{\Wtwo}{W_2}
\newcommand{\Wr}{W_r}
\newcommand{\dHel}{h}
\newcommand{\dV}{V}
\newcommand{\dK}{K}
\newcommand{\dPhi}{\rho_\phi}
\newcommand{\drPhi}{d_{\rho_\phi}}

\newcommand{\p}{p}
\newcommand{\hatp}{\hat{p}}
\newcommand{\pvec}{\mathbf{p}}
\newcommand{\qvec}{\mathbf{q}}
\newcommand{\vecp}{\mathbf{p}}
\newcommand{\vecq}{\mathbf{q}}
\newcommand{\Xcal}{\mathcal{X}}
\newcommand{\Pcal}{\mathcal{P}}
\renewcommand{\Qcal}{\mathcal{Q}}
\newcommand{\Qcall}{\mathcal{Q}(\mathbf{p},\mathbf{p}')}
\newcommand{\Ocal}{\mathcal{O}}
\newcommand{\Gcal}{\mathcal{G}}
\newcommand{\Gbar}{\bar{\mathcal{G}}}
\newcommand{\Diam}{\operatorname{Diam}}
\newcommand{\leqa}{\lesssim}
\newcommand{\geqa}{\gtrsim}
\newcommand{\supp}{\operatorname{supp}}

\newcommand{\tildeK}{\tilde{K}}
\newcommand{\tildef}{\tilde{f}}
\newcommand{\tildeg}{\tilde{g}}
\newcommand{\s}{s}
\newcommand{\kk}{t}

\newcommand{\m}{r}

\makeatother

\begin{document}
\begin{frontmatter}

\title{Convergence of latent mixing measures in finite and infinite
mixture models}
\runtitle{Convergence of mixing measures}

\begin{aug}
\author[A]{\fnms{XuanLong} \snm{Nguyen}\corref{}\thanksref{t1}\ead[label=e1]{xuanlong@umich.edu}}
\runauthor{X. Nguyen}
\affiliation{University of Michigan}
\address[A]{Department of Statistics\\
University of Michigan \\
456 West Hall \\
Ann Arbor, Michigan 48109-1107\\
USA \\
\printead{e1}} %adresu isvedimo komanda gale!
\end{aug}

\thankstext{t1}{Supported in part by NSF Grants CCF-1115769 and OCI-1047871.}

% HISTORY:
\received{\smonth{1} \syear{2012}}
\revised{\smonth{10} \syear{2012}}

% ABSTRACT
%
\begin{abstract}
This paper studies convergence behavior of latent mixing measures that arise
in finite and infinite mixture models, using transportation distances
(i.e., Wasserstein metrics).
The relationship between Wasserstein distances on the space of
mixing measures and $f$-divergence functionals such as Hellinger and
Kullback--Leibler distances on the space of mixture distributions
is investigated in detail using various identifiability conditions.
Convergence in Wasserstein metrics for discrete measures implies
convergence of individual atoms that provide support for the measures,
thereby providing a natural interpretation of convergence of clusters
in clustering applications where mixture models are typically employed.
Convergence rates
of posterior distributions for latent mixing measures are established,
for both finite mixtures of multivariate distributions and
infinite mixtures based on the Dirichlet process.
\end{abstract}

% KEYWORDS

\end{frontmatter}

%s1 #&#
\section{Introduction}
A notable feature in the development of hierarchical and Bayesian
nonparametric models is the role of mixing measures, which help to
combine relatively simple models into richer classes of statistical
models~\cite{Linsay-95,McLachlan-Basford-88}. In recent years the
mixture modeling methodology has been significantly extended by many
authors taking the mixing measure to be random and infinite-dimensional
via suitable priors constructed in a nested, hierarchical and
nonparametric manner. This results in rich models that can fit more
complex and high-dimensional data (see, e.g.,
\cite{Gelfand-etal-05,Teh-etal-06,Rodriguez-etal-08,Petrone-etal-09,Nguyen-10}
for several examples of such models, as well as a recent book
\cite{Hjort-etal-10}).

The focus of this paper is to analyze convergence behavior of the posterior
distribution of latent mixing measures as they arise in several
mixture models, including finite mixtures and the infinite Dirichlet
process mixtures.
%Although initially viewed as a modeling device for density estimation
%important role in the interpretation of the data population that they
%help
%to model, particularly in clustering applications.
Let $G = \sum_{i=1}^{k}p_i \delta_{\theta_i}$ denote a discrete probability
measure. Atoms $\theta_i$'s are elements in space $\Theta$,
while vector of probabilities $\pvec= (p_1,\ldots, p_k)$ lies in a
$k-1$-dimensional probability simplex. In a mixture setting, $G$ is
combined with a likelihood density $f(\cdot|\theta)$ with respect
to a dominating measure $\mu$ on $\Xcal$, to yield the mixture density:
$p_G(x) = \int f(x|\theta) \,dG(\theta)= \sum_{i=1}^{k}p_i
f(x|\theta_i)$.
In a clustering application, atoms $\theta_i$'s represent distinct
behaviors in a heterogeneous data population, while mixing
probabilities $p_i$'s are
the associated proportions of such behaviors.
Under this interpretation, there is a need for
comparing and assessing the quality of mixing measure $\hat{G}$
estimated on the basis of available data.
%A challenge in this regard,
%from both theoretical and methodological standpoints, appears
%to stem from the lack of suitable distance metrics for the space
%of discrete probability measures. This paper seeks to address this
%issue, by studying
%a class of distance functionals for discrete probability measures,
%and studying the convergence of the latent discrete mixing distribution
%in the proposed distance metrics.}
%
An important work in this direction is by Chen~\cite{Chen-95}, who
used the
$L_1$ metric on the cumulative distribution functions on the real line to
study convergence rates of the mixing measure $G$.
Chen's results were subsequently extended to a Bayesian
estimation setting for a univariate mixture
model~\cite{Ishwaran-James-Sun-01}. These works were limited to only
univariate and finite mixture models, with $k$ bounded by a known
constant, while our interest is when $k$ may be unbounded and
$\Theta$ is multidimensional or even an abstract space.
%For instance, $\Theta$ may be a subset of a function space
%as in the work of~\cite{Gelfand-etal-05,Nguyen-10}, or
%a space of probability measures~\cite{Rodriguez-etal-08}.

The analysis of consistency and convergence rates of posterior distributions
for Bayesian estimation has seen much progress in the past decade.
Key recent references
include \cite
{Barron-Shervish-Wasserman-99,Ghosal-Ghosh-vanderVaart-00,Shen-Wasserman-01,Walker-04,Ghosal-vanderVaart-07,Walker-Lijoi-Prunster-07}.
Analysis of specific mixture models in a Bayesian setting has also been
studied \cite
{Ghosal-Ghosh-Ramamoorthi-99,Genovese-Wasserman-00,Ishwaran-Zarepour-02,Ghosal-vanderVaart-07b}.
All these works primarily
focus on the convergence behavior of the posterior distribution
of the data density~$p_G$. On the other hand, results
concerned with the convergence behavior of latent mixing
measures $G$ are quite rare. Notably, the analysis
of convergence for mixing (smooth) densities often arises in the context
of frequentist estimation
for deconvolution problems, mainly within the
kernel density estimation method (e.g., \cite
{Carroll-Hall-88,Zhang-90,Fan-91}).
We also note recent progress on consistent parameter estimation
for certain finite mixture models, for example, in an overfitted
setting~\cite{Rousseau-Mengersen-11} or with an emphasis on computational
efficiency~\cite{Belkin-Sinha-10,Kalai-etal-12}.

The primary contribution of this paper is to show that the Wasserstein
distances provide a natural and useful metric for the analysis of
convergence for latent mixing measures in mixture models,
and to establish convergence rates of posterior distributions
in a number of well-known Bayesian nonparametric and mixture models.
Wasserstein distances originally arose in
the problem of optimal transportation~\cite{Villani-03}. Although not
as popular as well-known divergence functionals such as Kullback--Leibler,
total variation and Hellinger distances, Wasserstein distances have
been utilized in a number of statistical contexts
(e.g.,~\cite{Dudley-76,Mallows-72,Bickel-Freedman-81,delBarrio-etal-99}).
For discrete probability measures, they can be obtained by a
minimum matching (or moving) procedure between the
sets of atoms that provide support
for the measures under comparison, and consequentially are simple to
compute. Suppose that $\Theta$ is equipped with a metric $\rho$.
Let $G' = \sum_{j=1}^{k'}p'_j \delta_{\theta'_j}$. Then, for a given
$r \geq1$ the $L_r$
Wasserstein metric
on the space of discrete probability measures with support in $\Theta
$, namely,
$\Gbar(\Theta)$, is
\[
\Wr\bigl(G,G'\bigr) = \biggl[\inf_{\qvec} \sum
_{i,j} q_{ij}\rho^r\bigl(
\theta_i,\theta'_j\bigr)
\biggr]^{1/r},
\]
where the infimum is taken over all joint probability distributions on
$[1,\ldots,\break k]
\times[1,\ldots,k']$ such that $\sum_{j} q_{ij} = p_i$ and $\sum_{i} q_{ij} = p'_j$.

As clearly seen from this definition, Wasserstein distances inherit
directly the metric of
the space of atomic support $\Theta$, suggesting that they can be useful
for assessing estimation procedures for discrete measures
in hierarchical models. It is worth noting that if $(G_n)_{n\geq1}$ is a
sequence of discrete probability
measures with $k$ distinct atoms and $G_n$ tends to some discrete measure
$G_0$ in the $\Wr$
metric, then $G_n$'s ordered set of atoms must converge to $G_0$'s atoms
in $\rho$ after some permutation of atom labels.
Thus, in the clustering application illustrated above,
convergence of mixing measure $G$
may be interpreted as the convergence of distinct typical behavior
$\theta_i$'s
that characterize the heterogeneous data population.
A hint for the relevance of the Wasserstein distances
can be drawn from an observation that the $L_1$ distance for the CDFs
of univariate random variables, as studied by Chen~\cite{Chen-95}, is
in fact
a special case of the $W_1$ metric when $\Theta= \real$.

The plan for the paper is as follows.
Section~\ref{Sec-basic} investigates
the relationship between Wasserstein distances
for mixing measures and well-known divergence functionals for mixture densities
in a mixture model. We produce a simple lemma which gives
an upper bound on $f$-divergences between mixture densities by
certain Wasserstein distances between mixing measures. This
implies that $\Wr$ topology can be stronger than those
induced by divergences between mixture densities.
Next, we consider various identifiability
conditions under which convergence of mixture densities entails convergence
of mixing measures in a Wasserstein metric. We present two key
theorems, which provide upper bounds on $\Wtwo(G,G')$ in terms of
divergences between $p_G$ and $p_{G'}$.
Theorem~\ref{Thm-bound-supnorm} is applicable to
mixing measures with a bounded number of atomic support, generalizing
a result from~\cite{Chen-95}. Theorem~\ref{Thm-convolution}
is applicable to mixing measures with an unbounded number of support points,
but is restricted to only convolution mixture models.

Section~\ref{Sec-convergence} focuses on the convergence of posterior
distributions
of latent mixing measures in a Bayesian nonparametric setting.
Here, the mixing measure $G$ is endowed with a prior distribution $\Pi$.
Assuming an $n$-sample $X_1,\ldots,X_n$ that is generated according to
$p_{G_0}$,
we study conditions under which the postetrior distribution of $G$, namely,
$\Pi(\cdot|X_1,\ldots,X_n)$, contracts to the ``truth'' $G_0$ under the
$\Wtwo$ metric,
and provide the contraction rates.
In Theorems~\ref{Thm-convergence-1} and~\ref{Thm-convergence-2} of
Section~\ref{Sec-convergence}, we
establish the
convergence rates for the posterior\vadjust{\goodbreak} distribution for $G$ in terms of the
$\Wtwo$ metric.
These results are proved using the standard approach of Ghosal, Ghosh and
van der Vaart
\cite{Ghosal-Ghosh-vanderVaart-00}. Our convergence theorems
have several notable features. They rely on separate conditions
for the prior $\Pi$ and likelihood function $f$, which are typically
simpler to verify than
conditions formulated in terms of mixture densities.
The claim of convergence in Wasserstein metrics is typically
stronger than the weak convergence induced by the Hellinger metric
in the existing work mentioned above.

In Section~\ref{Sec-examples} posterior consistency and convergence rates
of latent mixing measures are derived, possibly for the first time,
for a number of well-known mixture
models in the literature, including finite mixtures of multivariate
distributions
and infinite mixtures based on Dirichlet processes.
For finite mixtures with a bounded number of atomic\vspace*{1pt} support in $\real^d$,
the posterior convergence rate for mixing measures is
$(\log n)^{1/4}n^{-1/4}$ under suitable identifiability conditions. This
rate is optimal up to a logarithmic factor in the minimax sense.
For Dirichlet process
mixtures defined on $\real^d$,
specific rates are established under smoothness conditions of the likelihood
density function~$f$. In particular, for ordinary smooth likelihood densities
with smoothness $\beta$ (e.g., Laplace), the rate achieved is
$(\log n/n)^\gamma$ for any $\gamma< \frac{2}{(d+2)(4+(2\beta+1)d)}$.
For supersmooth likelihood densities with smoothness $\beta$ (e.g., normal),
the rate achieved is $(\log n)^{-1/\beta}$.
%Finally, for finite mixtures of Gaussian processes, we are also able to
%establish a convergence rate by utilizing a result on
%(single) Gaussian process prior by van der Vaart and van Zanten

%pa1.subsection.subsubsection.1 #&#
\subsection*{Notation}
%measures. For discrete measures
%$G = \sum_{i=1}^{k}p_i\delta_{\theta_i}$ and $G'=
%the following mixture densities with respect to some
%dominating measure on $\Xcal$, respectively,
%$p_G(x) = \int f(x|\theta)\,dG(\theta) = \sum_{i=1}^{k}p_i f(x|
%$p_{G'}(x) = \int f(x|\theta')\,dG(\theta') = \sum_{j=1}^{k'}p'_j f(x|
%function on $\Xcal$.
%}
For ease of notation, we also use $f_i$ in place of $f(\cdot|\theta_i)$
and $f'_j$ in place of $f(\cdot|\theta'_j)$ for likelihood density functions.
Divergences (distances) studied in the paper include the total
variational distance:
$\dV(p_G,p_{G'}) = \frac{1}{2}\int|p_G(x) - p_{G'}(x)| \,d\mu(x)$,
Hellinger distance:
\[
\dHel^2(p_{G},p_{G'}) = \frac{1}{2}\int\bigl(\sqrt{p_G(x)} - \sqrt {p_{G'}(x)}\bigr)^2 \,d\mu(x)
\]
and Kullback--Leibler divergence:
\[
\dK(p_{G},p_{G'}) = \int p_{G}(x)\log\bigl(p_{G}(x)/p_{G'}(x)\bigr) \,d\mu(x).
\]
These divergences are related by $\dV^2/2
\leq\dHel^2 \leq\dV$ and $\dHel^2 \leq\dK/2$.
$N(\varepsilon,\Theta,\rho)$ denotes
the covering number of the metric space $(\Theta,\rho)$, that is, the
minimum number
of $\varepsilon$-balls needed to cover the entire space $\Theta$.
$D(\varepsilon,\Theta,\rho)$
denotes the packing number of $(\Theta,\rho)$, that is, the maximum
number of
points that are mutually separated by at least $\varepsilon$ in distance.
They are related by $N(\varepsilon,\Theta,\rho) \leq D(\varepsilon,\Theta
,\rho) \leq N(\varepsilon/2,\Theta,\rho)$. $\Diam(\Theta)$ denotes
the diameter of $\Theta$.

%s2 #&#
\section{Transportation distances for mixing measures}
\label{Sec-basic}

%s2.1 #&#
\subsection{Definition and a basic inequality}
Let $(\Theta,\rho)$ be a space equipped with a nonnegative distance function
$\rho\dvtx  \Theta\times\Theta\rightarrow\real_+$, that is, a function
that satisfies
$\rho(\theta_1,\theta_2) = 0 $ if and only if $\theta_1=\theta_2$.
If, in addition, $\rho$ is symmetric
($\rho(\theta_1,\theta_2) = \rho(\theta_2,\theta_1)$) and satisfies\vadjust{\goodbreak}
the triangle inequality, then it is a proper metric.
A discrete probability measure $G$ on a measure space equipped with the Borel
sigma algebra takes the form $G = \sum_{i=1}^{k} \p_i \delta_{\theta_i}$
for some $k \in\Nat\cup\{+\infty\}$, where $\mathbf{p} = (\p_1,\p_2,\ldots,
\p_k)$ denotes the proportion vector, while $\bolds{\theta} =
(\theta_1,\ldots,\theta_k)$
are the associated atoms in $\Theta$.
$\mathbf{p}$ has to satisfy $0\leq\p_i \leq1$ and $\sum_{i=1}^{k}\p_k = 1$.
[With a bit abuse of notation, we write $k=\infty$ when $G = \sum_{i=1}^{\infty}
p_i \delta_{\theta_i}$ has
countably infinite support points represented by the infinite
sequence of atoms $\bolds{\theta} = (\theta_1,\ldots)$ and the
associated sequence
of probability mass $\pvec$.]
Likewise, $G' = \sum_{j=1}^{k'}p'_j\delta_{\theta'_j}$ is another discrete
probability measure that has at most $k'$ distinct atoms.

Let $\Gcal_k(\Theta)$ denote the space of all discrete probability measures
with at most $k$ atoms. Let $\Gcal(\Theta) = \bigcup_{k\in\Nat_+}
\Gcal_k(\Theta)$,
the set of all discrete measures with finite support.
Finally, $\bar{\Gcal}(\Theta)$ denotes the space of all discrete measures
(including those with countably infinite support).

Let $\qvec=(q_{ij})_{i\leq k; j\leq k'} \in[0,1]^{k\times k'}$
denote a joint probability distribution on $\Nat_+\times\Nat_+$
that satisfies the marginal
constraints: $\sum_{i=1}^{k}q_{ij} = p'_j$
and $\sum_{j=1}^{k'}q_{ij} = p_i$ for any $i=1,\ldots,k; j=1,\ldots,k'$.
We also call $\qvec$ a coupling of $\pvec$ and $\mathbf{p}'$.
Let $\Qcall$ denote the space of all such couplings.
We start with the general transportation distance:
%
%de1 #&#
\begin{definition}
\label{Def-cd}
Let $\rho$ be a distance function on $\Theta$.
The transportation distance for two discrete
measures $G(\mathbf{p},\bolds{\theta})$ and $G'(\mathbf{p}',\bolds
{\theta}')$ is
%
%e1 #&#
\begin{equation}
d_\rho\bigl(G,G'\bigr) = \inf_{\qvec\in\Qcall} \sum
_{i,j} q_{ij}\rho \bigl(
\theta_i,\theta'_j\bigr).
\end{equation}
\end{definition}
When $\Theta$ is a metric space (e.g., $\real^d$)
and $\rho$ is taken to be its metric, we revert to the more standard notation
of Wasserstein metrics,
$W_1(G,G') \equiv d_\rho(G,G')$ and
$W_2^2(G,G') \equiv d_{\rho^2}(G,G')$. However,
$d_\rho$ will be employed when $\rho$ may be a general or a
nonstandard distance function or metric.

From here on, probability measure $G \in\Gbar(\Theta)$ plays
the role of the mixing distribution in a mixture model.
Let $f(x|\theta)$ denote the density (with respect to a dominating
measure $\mu$) of a random variable $X$ taking values in $\Xcal$,
given parameter $\theta\in\Theta$. For the ease of notation, we also use
$f_i(x)$ for $f(x|\theta_i)$.
Combining $G$ with the likelihood function $f$
yields a mixture distribution for $X$ that takes the following density:
\[
p_G(x) = \int f(x|\theta) \,dG(\theta)= \sum
_{i=1}^{k} \p_i f_i(x).
\]

A central theme in this paper is to explore
the relationship between Wasserstein distances
of mixing measures $G,G'$, for example, $d_\rho(G,G')$, and divergences
of mixture densities $p_G,p_{G'}$.
Divergences that play important roles in this paper\vadjust{\goodbreak}
are the total variational distance, the Hellinger
distance and the Kullback--Leibler distance.
All these are in fact instances of a broader class of divergences
known as the $f$-divergences (Csisz{\'a}r~\cite{Csi67}; Ali and Silvey~\cite{AliSil66}):
%
%de2 #&#
\begin{definition}
Let $\phi\dvtx  \real\rightarrow\real$ denote a convex function. An
$f$-divergence
(or Ali--Silvey distance) between two probability densities
$f_i$ and $f'_j$ is defined as $\dPhi(f_i,f'_j) = \int\phi(f'_j/f_i)
f_i \,d\mu$.
Likewise, the $f$-divergence between $p_G$ and $p_{G'}$ is
$\dPhi(p_G,p_{G'}) = \int\phi(p_{G'}/p_{G}) p_G \,d\mu$.
\end{definition}
$f$-divergences can be used as
a distance function or metric on $\Theta$.
When $\rho$ is taken to be an $f$-divergence, $\rho(\theta_i,\theta'_j)
:= \dPhi(f_i,f'_j)$, for a convex function $\phi$, we call
the corresponding transportation distance a
\textit{composite} transportation distance:
\[
d_{\rho_\phi}\bigl(G,G'\bigr):= \inf_{\qvec\in\Qcall}\sum
_{ij} q_{ij} \dPhi\bigl(f_i,f'_j
\bigr).
\]
For $\phi(u) = \frac{1}{2}(\sqrt{u}-1)^2$ we obtain the squared
Hellinger ($\rho_h^2 \equiv\dHel^2$), which induces the composite
transportation distance $d_{\rho_h^2}$.
For $\phi(u) = \frac{1}{2}|u-1|$ we obtain the variational
distance ($\rho_V \equiv\dV$), which induces $d_{\rho_V}$.
For $\phi(u) = -\log u$, we obtain
the Kullback--Leiber divergence ($\rho_K \equiv\dK$), which
induces~$d_{\rho_K}$.
\begin{lemma}
\label{Lem-gen-ineq}
Let $G,G' \in\Gbar(\Theta)$ such that both $\dPhi(p_G,p_{G'})$ and
$\drPhi(G,G')$ are finite for some convex function $\phi$.
Then, $\dPhi(p_G,p_{G'}) \leq\drPhi(G,G')$.
%In particular, $d_{h}^2(p_G,p_{G'}) \leq d_{\rho h^2}(G,G')$,
%$d_{V}(p_G,p_{G'}) \leq d_{\rho V} (G,G')$, and
%$d_{K}(p_G,p_{G'}) \leq d_{\rho K}(G,G')$.
\end{lemma}

%%%%%%%%%%%%%%%%%%%%%%%%%%%%%%%%%%%%%%%%%%%%%%%%%%%%%%%%%%%%%%%%%%%%%%%
%Collecting the results from the previous two lemmas
%we obtain the following inequalities for
%various distances between discrete measures $G$ and $G'$:
%%
%0\leq\dV^2(p_G,p_{G'})/2 \leq\dHel^2 (p_G,p_{G'}) \leq\{
%
%and $\{\dHel^2(p_G,p_{G'}), d_{\rho h}^2(G,G')\} \leq d_{\rho
%h^2}(G,G')\leq d_{\rho h}(G,G') \leq1$.
%
%$(G_n,G'_n)$.
%If $d_{\rho h}(G_n,G'_n) \rightarrow0$, then $d_{h}(p_{G_n},p_{G'_n})
%and $\dV(p_{G_n},p_{G'_n}) \rightarrow0$. In addition, under
%condition A,
%if $d_{\rho h}(G_n,G'_n) \rightarrow1$, then $d_{h}(p_{G_n},p_{G'_n})
%}
%%%%%%%%%%%%%%%%%%%%%%%%%%%%%%%%%%%%%%%%%%%%%%%%%%%%%%%%%%%%%%%%%%%%%%%
This lemma highlights a simple direction in the aforementioned
relationship: any $f$-divergence between mixture distributions $p_G$
and $p_{G'}$ is
dominated by a transportation distance between mixing measures $G$ and $G'$.
As will be evident in the sequel,
this basic inequality is also handy in enabling us to obtain upper bounds
on the power of tests. It also proves useful for establishing
lower bounds on small Kullback--Leibler ball probabilities in the space
of mixture
densities $p_G$ in terms of small ball probabilities in the metric
space $(\Theta,\rho)$.
The latter quantities are typically easier to obtain estimates for than
the former.
%
%ex1 #&#
\begin{example}
\label{example-gaussian}
Suppose that $\Theta= \real^d$, $\rho$ is the Euclidean metric,
$f(x|\theta)$ is the multivariate normal density $N(\theta,I_{d\times
d})$ with mean
$\theta$ and identity covariance matrix, then
$\dHel^2(f_i,f'_j) = 1- \exp-\frac{1}{8}\|\theta_i-\theta'_j\|^2
\leq\frac{1}{8}\|\theta_i-\theta'_j\|^2 = \rho^2(\theta_i,\theta'_j)^2/8$. So,
$d_{\rho_h^2}(G,G') \leq d_{\rho^2}(G,G')/8$. The above lemma then
entails that
$\dHel^2(p_G,p_{G'}) \leq d_{\rho^2}(G,G')/8 = \Wtwo^2(G,G')/8$.

Similarly, for the Kullback--Leibler divergence, since
$\dK(f_i,f'_j) = \frac{1}{2}\|\theta_i-\theta'_j\|^2$,
by Lemma~\ref{Lem-gen-ineq},
$\dK(p_{G},p_{G'}) \leq d_{\rho_K}(G,G') = \frac{1}{2}d_{\rho^2}(G,G')
= \Wtwo(G,G')^2/2$.
%Next, suppose that $\Theta$ is a compact subset of $\real^d$
%and consider $\phi(u) = (\log u)^2$, which is a convex function
%on $[0,\infty)$. We have $\int f_i (\log f_i/f'_j)^2
%= O(\|\theta_i-\theta'_j)\|^2$, so $\int p_{G}[\log(p_G/p_{G'})]^2

For another example, if $f(x|\theta)$ is a Gamma density with location
parameter~$\theta$, $\Theta$ is a compact subset of $\real$ that is
bounded away from 0.
Then $\dK(f_i,f'_j) = O(|\theta_i-\theta_j|)$. This entails that
$\dK(p_G,p_{G'}) \leq d_{\rho_K}(G,G') \leq O(W_1(G,G'))$.
\end{example}

\subsection{Wasserstein metric identifiability in finite mixture models}
\label{Sec-identifiability}

Lemma~\ref{Lem-gen-ineq} shows that for many choices of $\rho$,
$d_\rho$ yields a stronger topology on $\Gbar(\Theta)$
than the topology induced by $f$-divergences
on the space of mixture distributions $p_G$. In other words,
convergence of $p_G$ may not imply convergence of $G$ in
transportation distances.
To ensure this property, additional conditions are
needed on the space of probability measures $\Gbar(\Theta)$, along
with identifiability
conditions for the family of likelihood functions $\{f(\cdot|\theta),
\theta\in\Theta\}$.

The classical definition of Teicher~\cite{Teicher-61} specifies
the family $\{f(\cdot|\theta), \theta\in\Theta\}$ to be identifiable
if for any $G,G'\in\Gcal(\Theta)$, $\|p_{G}-p_{G'}\|_\infty= 0$
implies that $G=G'$. We need a slightly stronger version,
allowing for the inclusion for discrete measures with infinite support:
%
%de3 #&#
\begin{definition}
The family $\{f(\cdot|\theta), \theta\in\Theta\}$ is finitely identifiable
if for any $G\in\Gcal_\Theta$ and $G'\in\bar{\Gcal}_\Theta$,
$|p_G(x)-p_{G'}(x)| = 0$ for almost all $x\in\Xcal$ implies that $G=G'$.
\end{definition}

To obtain convergence rates, we also need the notion of strong
identifiability of~\cite{Chen-95}, herein adapted to a multivariate setting.
%
%de4 #&#
\begin{definition}
\label{Def-identifiability}
Assume that $\Theta\subseteq\real^d$ and $\rho$ is the Euclidean metric.
The family $\{f(\cdot|\theta), \theta\in\Theta\}$
is strongly identifiable if $f(x|\theta)$ is twice differentiable
in $\theta$ and for any finite $k$ and $k$ different $\theta_1,\ldots,\theta_k$,
the equality
%
%e2 #&#
\begin{equation}
\label{Eqn-uniform} \ess\sup_{x\in\Xcal} \Biggl|\sum_{i=1}^{k}
\alpha_i f(x|\theta_i) + \beta_i^T
Df(x|\theta_i) + \gamma_i^T
D^2f(x|\theta_i) \gamma_i \Biggr| = 0
\end{equation}
implies that $\alpha_i = 0$, $\beta_i = \gamma_i = \mathbf{0} \in
\real^d$ for $i=1,\ldots,k$. Here, for each $x$,
$Df(x|\theta_i)$ and $D^{2}f(x|\theta_i)$
denote the gradient and the Hessian at $\theta_i$
of function $f(x|\cdot)$, respectively.
\end{definition}

Finite identifiability is satisfied for the family of Gaussian
distributions for both mean and variance parameters
\cite{Teicher-60}; see also Theorem 1 of~\cite{Ishwaran-Zarepour-02}.
Chen identified a broad class of families, including the Gaussian
family, for which the strong identifiability condition holds
\cite{Chen-95}.

Define $\psi(G,G') = \sup_{x}|p_G(x)-p_{G'}(x)|/\Wtwo^2(G,G')$ if
$G\neq G'$
and $\infty$ otherwise. Also define
$\psi_1(G,G') = \dV(p_G,p_{G'})/\Wtwo^2(G,G')$ if $G\neq G'$
and $\infty$ otherwise. The notion of strong identifiability is useful
via the following key result, which generalizes Chen's result to
$\Theta$ of arbitrary dimensions.\vfill\eject
%
%th1 #&#
\begin{theorem}[(Strong identifiability)]
\label{Thm-bound-supnorm}
Suppose that $\Theta$ is a compact subset of~$\real^d$,
the family $\{f(\cdot|\theta),\theta\in\Theta\}$
is strongly identifiable, and for all $x \in\Xcal$, the Hessian
matrix $D^{2}f(x|\theta)$ satisfies a uniform Lipschitz condition
%
%e3 #&#
\begin{equation}
\label{Eqn-Lipschitz} \bigl|\gamma^T\bigl(D^{2}f(x|
\theta_1) - D^{2}f(x|\theta_2)\bigr)\gamma\bigr|
\leq C\|\theta_1-\theta_2\|^{\delta}\|\gamma
\|^2
\end{equation}
for all $x, \theta_1,\theta_2$ and some fixed $C$ and $\delta> 0$. Then,
for fixed $G_0 \in\Gcal_k(\Theta)$, where $k<\infty$,
%
%e4 #&#
\begin{equation}
\label{Eqn-bound-supnorm} \lim_{\varepsilon\rightarrow0} \inf_{G,G'\in\Gcal_k(\Theta)} \bigl\{ \psi
\bigl(G,G'\bigr)\dvtx  \Wtwo(G_0,G) \vee \Wtwo
\bigl(G_0,G'\bigr) \leq\varepsilon \bigr\} > 0.
\end{equation}
The assertion also holds with $\psi$ being replaced by $\psi_1$.
\end{theorem}

%pa2.2.subsubsection.1 #&#
\begin{Remark*}
Suppose that $G_0$ has exactly $k$ distinct support points in $\Theta$
(i.e., $G = \sum_{i=1}^{k}p_i \delta_{\theta_i}$ where
$p_i>0$ for all $i=1,\ldots,k$). Then, an examination of
the proof reveals that the requirement that $\Theta$ be compact is not needed.
Indeed, if there is a sequence of $G_n \in\Gcal_k(\Theta)$ such that
$W_2(G_0,G_n) \rightarrow0$,
then it is simple to show that
%using an argument in the first paragraph of the proof of Lemma
%
there is a subsequence of $G_n$ that also has $k$ distinct atoms, which
converge in the $\rho$ metric to the set of $k$ atoms of $G_0$
(up to some permutation of the labels). The proof of the theorem
proceeds as before.

%(ii) In Section~\ref{normspace} we briefly discuss the extension of
%the notion of strong identifiability to an infinite dimensional
%setting via first and
%second order Fr\'echet derivatives in normed spaces.

For the rest of this paper, by strong identifiability we always mean conditions
specified in Theorem~\ref{Thm-bound-supnorm}
so that equation (\ref{Eqn-bound-supnorm}) can be deduced. This
practically means that the conditions
specified by (\ref{Eqn-uniform}) and (\ref
{Eqn-Lipschitz}) be given, while
the compactness of $\Theta$ may sometimes be required.
\end{Remark*}

%(ii) When $f(x|\theta)$ is a multivariate normal
%density (as in Example 5) we also know that $\Wtwo^2(G,G')
%h^2}(G,G')$ for
%some small constant $c> 0$.
%This clarifies the relationship between $\dV(p_G,p_{G'})$
%and $d_{\rho h^2}(G,G')$ reported empirically in Fig
%}

%s2.3 #&#
\subsection{Wasserstein metric identifiability in infinite mixture models}
\label{Sec-identifiability-2}

Next, we state a counterpart of Theorem~\ref{Thm-bound-supnorm} for
$G,G'\in\Gbar(\Theta)$, that is, mixing measures with a potentially
unbounded number of support points. We restrict our attention to
convolution mixture models on $\real^d$. That is, the likelihood
density function $f(x|\theta)$, with respect to Lebesgue, takes the
form $f(x-\theta)$ for some multivariate density function $f$ on
$\real^d$. Thus, $p_G(x) = G*f(x) = \sum_{i=1}^{k} p_i f(x-\theta_i)$
and $p_{G'}(x) = G'*f(x) = \sum_{j=1}^{k'} p'_j
f(x-\theta'_{j})$.\vspace*{1pt}

The key assumption is concerned with the smoothness of density function
$f$. This is
characterized in terms of the tail behavior of the Fourier transform
$\tilde{f}$
of $f\dvtx  \tilde{f}(\omega) = \int_{\real^d} e^{-i\langle\omega,
x\rangle} f(x) \,dx$.
We consider both ordinary smooth densities (e.g., Laplace and Gamma) and
supersmooth densities (e.g., normal).

%th2 #&#
\begin{theorem}
\label{Thm-convolution}
Suppose that $G,G'$ are probability measures that place full support on
a bounded subset $\Theta\subset\real^d$.
$f$ is a density function on $\real^d$ that is symmetric
(around 0), that is, $\int_{A} f \,dx = \int_{-A} f \,dx$ for any Borel
set $A\subset\real^d$.
Moreover, assume that $\tildef(\omega) \neq0$ for all $\omega\in
\real^d$.
\begin{longlist}[(2)]
\item[(1)] \textup{Ordinary smooth likelihood.}
Suppose that\vspace*{1pt} $|\tildef(\omega) \prod_{j=1}^{d}|\omega_j|^\beta|
\geq d_0$ as $\omega_j \rightarrow\infty$ $(j=1,\ldots, d)$ for
some positive constants $d_0$ and $\beta$. Then for any
$m< 4/(4+(2\beta+1)d)$, there\vadjust{\goodbreak} is some constant $C(d,\beta, m)$
dependent only
on $d, \beta$ and $m$ such that
\[
\Wtwo^2\bigl(G,G'\bigr) \leq C(d,\beta, m)
\dV(p_G,p_{G'})^m
\]
as $\dV(p_{G},p_{G'}) \rightarrow0$.

\item[(2)] \textup{Supersmooth likelihood.}
Suppose that $|\tildef(\omega) \prod_{j=1}^{d}\exp(|\omega_j|^{\beta}/\gamma)| \geq d_0$
as $\omega_j \rightarrow\infty$ $(j=1,\ldots, d)$
for some positive constants $\beta,\gamma,d_0$. Then there is some
constant $C(d,\beta)$
dependent only on $d$ and $\beta$ such that
\[
\Wtwo^2\bigl(G,G'\bigr) \leq C(d,\beta) \bigl(-\log
\dV(p_G,p_{G'})\bigr)^{-2/\beta}
\]
as $\dV(p_{G},p_{G'}) \rightarrow0$.
\end{longlist}
\end{theorem}
%
%pa2.3.subsubsection.1 #&#
\begin{Remark*}
The theorem does not actually require
that mixing measures $G,G'$ be discrete. Moreover, from the proof of
the theorem, the
condition that the support points of $G$ and $G'$ lie in a bounded
subset of $\real^d$ can be removed and replaced by
the boundedness of a given moment of the mixing measures. The upper bound
remains the same for the supersmooth likelihood case. For the ordinary smooth
case, we obtain a slightly weaker upper bound for $W_2(G,G')$.
\end{Remark*}
%
%ex2 #&#
\begin{example} For the standard normal density on $\real^d$, $\tildef
(\omega)
= \prod_{j=1}^{d}\exp-\break\omega_i^2/2$, we obtain that
$\Wtwo^2(G,G') \leqa(-\log\dV(p_G,p_{G'}))^{-1}$ as
$\Wtwo(G,G') \rightarrow0$ [so that $\dV(p_G,p_{G'}) \rightarrow0$, by
Lemma~\ref{Lem-gen-ineq}]. For a\vspace*{1pt} Laplace density on $\real$, for example,
$\tildef(\omega) = \frac{1}{1+\omega^2}$, then $\Wtwo^2(G,G')
\leqa\dV(p_G,p_{G'})^m$ for any $m < 4/9$, as $\Wtwo(G,G')
\rightarrow0$.
\end{example}

%s3 #&#
\section{Convergence of posterior distributions of mixing measures}
\label{Sec-convergence}

We turn to a study of convergence of mixing measures
in a Bayesian setting. Let $X_1,\ldots,X_n$ be an i.i.d. sample according
to the
mixture density $p_G(x) = \int f(x|\theta) \,dG(\theta)$, where $f$ is
known, while
$G= G_0$ for some unknown mixing measure in $\Gcal_k(\Theta)$. The
true number of
support points for $G$ may be unknown (and/or unbounded). In the
Bayesian estimation
framework, $G$ is endowed with
a prior distribution $\Pi$ on a suitable measure space of discrete
probability measures in $\Gbar(\Theta)$.
The posterior distribution of $G$ is given by, for any measurable set
$B$,
\[
\Pi(B|X_1,\ldots,X_n) = \int_{B}
\prod_{i=1}^{n}p_G(X_i)
\,d\Pi(G) \Big/ \int\prod_{i=1}^{n}p_G(X_i)
\,d\Pi(G).
\]

We shall study conditions under which the posterior distribution is
consistent, that is,
it concentrates on arbitrarily small $\Wtwo$ neighborhoods of $G_0$,
and establish the rates of the convergence.
We follow the general framework of Ghosal, Ghosh and
van der Vaart~\cite{Ghosal-Ghosh-vanderVaart-00},
who analyzed convergence behavior of posterior distributions
in terms of $f$-divergences such as Hellinger and variational distances
on the mixture densities of the data. In the following
we formulate two convergence theorems\vadjust{\goodbreak} for the mixture model setting
(which can be viewed as counterparts of Theorems 2.1 and 2.4
of~\cite{Ghosal-Ghosh-vanderVaart-00}). A notable feature
of our theorems is that conditions (e.g., entropy and prior concentration)
are stated directly in terms of the Wasserstein metric, as opposed to
$f$-divergences on the mixture densities. They may be typically
separated into
independent conditions for the prior for $G$ and the likelihood family
and are
simpler to verify for mixture models.

The following notion plays a central role in our general results.
%
%de5 #&#
\begin{definition} Fix $G_0 \in\Gbar(\Theta)$.
Let $\Gcal\subset\Gbar(\Theta)$. Define the Hellinger information
of the $\Wtwo$ metric for subset $\Gcal$ as a real-valued function on the
real line $\Psi_\Gcal\dvtx \real\rightarrow\real$:
%
%e5 #&#
\begin{equation}
\label{Eqn-def-Ck} \Psi_\Gcal(r) = \inf_{G\in\Gcal: \Wtwo(G_0,G) \geq r/2}
\dHel^2(p_{G_0},p_{G}).
\end{equation}
\end{definition}
Note the dependence of $\Psi_\Gcal$ on the (fixed) $G_0$, but this is
suppressed
for ease of notation.
It is obvious that $\Psi_\Gcal$ is a nonnegative and nondecreasing
function. The following
characterizations of $\Psi_\Gcal$ are simple consequences of Theorems
\ref{Thm-bound-supnorm}
and~\ref{Thm-convolution}:
%
%pr1 #&#
\begin{proposition}
\label{Thm-1}
\textup{(a)} Suppose that $G_0 \in\Gcal_k(\Theta)$, and both $\Gcal_k(\Theta)$
and $\Gcal$ are compact in the Wasserstein
topology. In addition, assume that the family of likelihood functions
is finitely identifiable.
Then, $\Psi_\Gcal(r) > 0$ for all $r>0$.

\textup{(b)} Suppose that $\Theta\subset\real^d$ is compact, and the
family of likelihood
functions is strongly identifiable as specified in Theorem \ref
{Thm-bound-supnorm}.
Then, for \textup{each} $k$ there is a constant $c(k,G_0) >0$ such that
$\Psi_{\Gcal_k(\Theta)}(r) \geq c(k,G_0)r^4$ for all $r > 0$.
%Furthermore, the compactness of $\Theta$
%can be replaced by the assumption that $G_0$ has exactly $k$ distinct
%atoms, and that
%$\Theta$ is bounded.
%functions is strongly identifiable,
%and given a fixed $G_0$ that has exactly $k$ distinct atoms in $

\textup{(c)} Suppose that $\Theta\subset\real^d$ is bounded,
and the family of likelihood functions is
ordinary smooth with parameter $\beta$, as specified in Theorem \ref
{Thm-convolution}. Then,\vspace*{1pt} for any $d'>d$ there is some constant
$c(d,\beta)$ such that\vspace*{1pt}
$\Psi_{\Gbar(\Theta)}(r) \geq c(d,\beta) r^{4+(2\beta+1)d'}$
for all \mbox{$r>0$}.
For the supersmooth likelihood family, we have $\Psi_{\Gbar(\Theta
)}(r) \geq\exp[-c(d,\beta)r^{-\beta}]$
for all $r>0$.
\end{proposition}

A main ingredient in the analysis of convergence of posterior distributions
is through proving the existence of tests for subsets of parameters of interest.
A test $\varphi_n$ is a measurable indicator function of the i.i.d.
sample $X_1,\ldots,X_n$.
For a fixed pair of measures $(G_0,G_1)$ such that $G_1 \in\Gcal$,
where $\Gcal$ is a given subset of $\Gbar(\Theta)$,
consider tests for discriminating $G_0$ against a
closed Wasserstein ball centered at~$G_1$. Write
\[
B_W(G_1,r) = \bigl\{G \in\Gbar(\Theta)\dvtx
\Wtwo(G_1,G) \leq r\bigr\}.
\]
The following lemma highlights the role of
the Hellinger information:
%
%le2 #&#
\begin{lemma}
\label{Lem-convexball}
Suppose that $(\Theta,\rho)$ is a metric space. Fix $G_0 \in\Gbar
(\Theta)$
and consider $\Gcal\subset\Gbar(\Theta)$.
Given $G_1 \in\Gcal$, let $r= \Wtwo(G_0,G_1)$.
Suppose that either one of the following two sets of conditions holds:\vadjust{\goodbreak}
\begin{longlist}[(II)]
\item[(I)] $\Gcal$ is a convex set, in which case,
let $M(\Gcal,G_1,r) = 1$.

\item[(II)] $\Gcal$ is nonconvex, while $\Theta$ is a totally bounded
and bounded set. In addition, for some constants
$C_1 >0,\alpha\geq1$, $\dHel(f_i,f'_j) \leq C_1 \rho^\alpha(\theta_i,\theta'_j)$ for
any likelihood functions $f_i,f'_j$ in the family. In this case, define
%
%e6 #&#
\begin{equation}
\label{Eqn-M} M(\Gcal,G_1,r) = D \biggl(\frac{\Psi_\Gcal(r)^{1/2}}{2\Diam
(\Theta)^{\alpha-1}\sqrt{C_1}}, \Gcal
\cap B_W(G_1,r/2), \Wtwo \biggr).
\end{equation}
\end{longlist}
Then, there exist tests $\{\varphi_n\}$ that have the following properties:
%
%e7 #&#
%e8 #&#
\begin{eqnarray}
\label{Eqn-type-1} P_{G_0} \varphi_n & \leq& M(
\Gcal,G_1,r) \exp\bigl[-n \Psi_\Gcal (r)/8\bigr],
\\
\label{Eqn-type-2}
\sup_{G\in\Gcal\cap B_W(G_1,r/2)} P_{G}(1-\varphi_n) & \leq& \exp
\bigl[-n \Psi_\Gcal(r)/8\bigr].
\end{eqnarray}
Here, $P_{G}$ denotes the expectation under the mixture distribution
given by density~$p_{G}$.
\end{lemma}

%pa3.subsection.subsubsection.1 #&#
\begin{Remark*}
The set of conditions (II) is needed when $\Gcal$ is not convex, an example
of which is $\Gcal= \Gcal_k(\Theta)$, the space of measures with at most
$k$ support points in $\Theta$. It is interesting to note that the
loss in
test power due to the lack of convexity is captured by the local
entropy term $\log M(\Gcal,G_1,r)$. This quantity is defined in terms of
the packing by small $\Wtwo$ balls whose radii are specified by
the Hellinger information. Hence, this packing number
can be upper bounded by exploiting a lower bound of the
Hellinger information.
%We note briefly that this quantity can be
%easily bounded in some cases: for instance,
%if $\Theta\subset\real^d$ and $\Gcal= \Gcal_k(\Theta)$, then by
%Proposition
%$\Psi_\Gcal(r) \geq c(k,G_0)r^4$.
%In addition, under
%the assumption that for all $G\in\Gcal$, both the masses for
%all $k$ support points of $G$, and the pairwise distances of such
%support points,
%are uniformly bounded away from 0,
%it can be shown that $M(\Gcal,G_1,r)$ is bounded by a constant that
%depends on only
%$G_0, \Theta$ and $k$ (Further details are elaborated in the proof of
%Theorem
\end{Remark*}

Next, the existence of the test can be shown for discriminating $G_0$
against the
complement of a closed Wasserstein ball:
%(see also
%
%le3 #&#
\begin{lemma}
\label{Lem-ballcomplement}
Assume that conditions of Lemma~\ref{Lem-convexball} hold and define
$M(\Gcal,\break G_1,r)$ as in Lemma~\ref{Lem-convexball}.
Suppose that for some nonincreasing function $D(\varepsilon)$, some
$\varepsilon_n \geq0$ and every $\varepsilon> \varepsilon_n$,
%
%e9 #&#
\begin{equation}
\label{Eqn-entropy-1}\quad \sup_{G_1\in\Gcal}M(\Gcal,G_1,\varepsilon)
\times D\bigl(\varepsilon/2, \Gcal\cap B_W(G_0,2\varepsilon)
\setminus B_W(G_0,\varepsilon), %\{G \in\Gcal: \varepsilon\leq d_{\rho}(G_0,G) \leq2\varepsilon\},
\Wtwo\bigr)
\leq D(\varepsilon).
\end{equation}
Then, for every $\varepsilon> \varepsilon_n$, for any $t_0\in\Nat$,
there exist tests $\varphi_n$
(depending on $\varepsilon> 0$) such that
%
%e10 #&#
%e11 #&#
\begin{eqnarray}
\label{Eqn-type-1b} P_{G_0} \varphi_n & \leq& D(\varepsilon)
\sum_{t=t_0}^{\lceil\Diam(\Theta)/\varepsilon\rceil}\exp\bigl[-n
\Psi_{\Gcal}(t\varepsilon)/8\bigr],
\\
\label{Eqn-type-2b}\quad \sup_{G\in\Gcal\dvtx  \Wtwo(G_0,G) > t_0\varepsilon} P_{G} (1-
\varphi_n) & \leq& \exp\bigl[-n \Psi_{\Gcal}(t_0
\varepsilon)/8\bigr].
\end{eqnarray}
\end{lemma}

%pa3.subsection.subsubsection.2 #&#
\begin{Remark*}
It is interesting to observe that function $D(\varepsilon)$ is used to
control the
packing number of thin layers of Wasserstein balls (a similar quantity that
also arises via the peeling argument in \cite
{Ghosal-Ghosh-vanderVaart-00} (Theorem 7.1)),
in addition to the packing number $M$ of small\vadjust{\goodbreak} Wasserstein balls in
terms of smaller Wasserstein balls whose radii are specified
in terms of the Hellinger information function. As in the previous
lemma, the
latter packing number appears intrinsic to the analysis of convergence
for mixing measures.

The preceeding two lemmas provide the core argument for establishing
the following general posterior contraction theorems for latent
mixing measures in a mixture model.
The following two theorems have three types of conditions. The first
is concerned with the size of support of $\Pi$, often quantified
in terms of its entropy number. Estimates of the entropy number defined
in terms of Wasserstein metrics for several measure classes of interest
are given in Lemma~\ref{Lem-entropy}.
The second condition is on the Kullback--Leibler support
of $\Pi$, which is related to both the space of discrete measures
$\Gbar(\Theta)$
and the family of likelihood functions
$f(x|\theta)$. The Kullback--Leibler neighborhood is defined as
%
%e12 #&#
\begin{equation}
\label{Eqn-kl-neighborhood}\quad B_K(\varepsilon) = \biggl\{G \in\Gbar(
\Theta)\dvtx  -P_{G_0} \biggl(\log \frac{p_{G}}{p_{G_0}} \biggr) \leq
\varepsilon^2, P_{G_0} \biggl( \log\frac{p_{G}}{p_{G_0}}
\biggr)^2 \leq\varepsilon^2 \biggr\}.
\end{equation}
The third type of condition is on the Hellinger information of the
$\Wtwo$ metric,
function $\Psi_\Gcal(r)$, a characterization of which is given above.
\end{Remark*}
%
%th3 #&#
\begin{theorem}
\label{Thm-convergence-1}
Fix $G_0 \in\Gbar(\Theta)$, and assume that
the family of likelihood functions is finitely identifiable.
Suppose that for a sequence $(\varepsilon_n)_{n\geq1}$ that tends to a
constant (or 0)
such that $n\varepsilon_n^2 \rightarrow\infty$, and a constant $C>0$,
and \textup{convex} sets $\Gcal_n \subset\Gbar(\Theta)$, we have
%
%e13 #&#
%e14 #&#
%e15 #&#
\begin{eqnarray}
\label{Eqn-entropy-2s} \log D(\varepsilon_n, \Gcal_n,
\Wtwo) &\leq& n\varepsilon_n^2,
\\
\label{Eqn-support-1s} \Pi\bigl(\Gbar(\Theta)\setminus\Gcal_n\bigr)
&\leq&\exp\bigl[-n\varepsilon_n^2(C+4)\bigr],
\\
\label{Eqn-support-2s} \Pi\bigl(B_K(\varepsilon_n)\bigr)
&\geq&\exp\bigl(-n\varepsilon_n^2C\bigr).
\end{eqnarray}
Moreover, suppose $M_n$ is a sequence such that
%
%e16 #&#
%e17 #&#
\begin{eqnarray}
\label{Eqn-rate-conds-2s} \Psi_{\Gcal_n}(M_n
\varepsilon_n) &\geq&8\varepsilon_n^2(C+4),
\\
\label{Eqn-hellinger-sum} \exp\bigl(2n\varepsilon_n^2
\bigr)\sum_{j\geq M_n} \exp\bigl[-n\Psi_{\Gcal
_n}(j
\varepsilon_n)/8\bigr] &\rightarrow&0.
\end{eqnarray}
Then, $\Pi(G\dvtx  \Wtwo(G_0,G) \geq M_n\varepsilon_n|X_1,\ldots,X_n)
\rightarrow0$ in $P_{G_0}$-probability.
\end{theorem}

The following theorem uses a weaker condition on the
covering number, but it contains an additional condition on
the likelihood functions which may be
useful for handling the case of nonconvex sieves $\Gcal_n$.
%
%th4 #&#
\begin{theorem}
\label{Thm-convergence-2}
Fix $G_0 \in\Gbar(\Theta)$. Assume the following:
\begin{longlist}[(a)]
\item[(a)] The family of likelihood functions is finitely identifiable
and satisfies
$h(f_i,f'_j) \leq C_1 \rho^\alpha(\theta_i,\theta'_j)$ for
any\vspace*{1pt}
likelihood functions
$f_i,f'_j$ in the family, for some constants $C_1 >0, \alpha\geq1$.
\item[(b)] There is a sequence of sets $\Gcal_n \subset\Gbar(\Theta)$
for which $M(\Gcal_n,G_1,\varepsilon)$ is defined by (\ref{Eqn-M}).
\item[(c)] There is a sequence $\varepsilon_n \rightarrow0$
such that $n\varepsilon_n^2$ is bounded away from 0 or tending to infinity,
and a sequence $M_n$ such that
%
%e18 #&#
%e19 #&#
%e20 #&#
%e21 #&#
\begin{eqnarray}\label{Eqn-entropy-2}
&& \log D\bigl(\varepsilon/2, \Gcal_n \cap B_W(G_0,2
\varepsilon) \setminus B_W(G_0,\varepsilon), \Wtwo\bigr)
\nonumber\\[-8pt]\\[-8pt]
&&\qquad{} +\sup_{G_1\in\Gcal_n} \log M(\Gcal_n,
G_1,\varepsilon) \leq n\varepsilon_n^2\qquad
\forall\varepsilon\geq\varepsilon_n,\nonumber\\[-12pt]\nonumber
\end{eqnarray}
\begin{eqnarray}
\label{Eqn-support-1}
\frac{\Pi(\Gbar(\Theta)\setminus\Gcal_n)}{\Pi(B_K(\varepsilon_n))} &=&
o\bigl(\exp\bigl(-2 n
\varepsilon_n^2\bigr)\bigr),
\\
\label{Eqn-support-2}\qquad
\frac{\Pi( B_W(G_0,2j\varepsilon_n) \setminus B_W(G_0, j\varepsilon_n))} {
\Pi(B_K(\varepsilon_n))} &\leq&\exp\bigl[n
\Psi_{\Gcal_n}(j\varepsilon_n)/16\bigr] \qquad\forall j\geq
M_n,\\[-12pt]\nonumber
\end{eqnarray}
\begin{eqnarray}
\label{Eqn-rate-cond}
\exp\bigl(2n\varepsilon_n^2\bigr)
\sum_{j\geq M_n} \exp\bigl[-n\Psi_{\Gcal
_n}(j
\varepsilon_n)/16\bigr] &\rightarrow&0.
\end{eqnarray}
\end{longlist}

Then, we have that
$\Pi(G\dvtx  \Wtwo(G_0,G) \geq M_n\varepsilon_n|X_1,\ldots,X_n)
\rightarrow0$ in $P_{G_0}$-probability.
\end{theorem}

%pa3.subsection.subsubsection.3 #&#
\begin{Remark*}
(i) From the theorem's proof, the above statement continues to
hold if conditions
(\ref{Eqn-support-2}) and (\ref{Eqn-rate-cond}) are replaced by the following
condition:
%
%e22 #&#
\begin{equation}
\exp\bigl(2n\varepsilon_n^2\bigr)/\Pi\bigl(B_K(
\varepsilon_n)\bigr) \sum_{j\geq M_n}\exp
\bigl[-n\Psi_{\Gcal_n}(j\varepsilon_n)/16\bigr] \rightarrow0.
\end{equation}

(ii) In Theorem~\ref{Thm-convergence-2} and in Theorem
\ref{Thm-convergence-1} augmented with condition (a) of
Theorem~\ref{Thm-convergence-2}, it is simple to deduce the posterior
convergence rates for the mixture density $p_G$. We obtain that
for a sufficiently large constant $M > 0$,
\[
\Pi\bigl(G\dvtx  h(p_{G_0},p_{G}) \geq M\varepsilon_n
| X_1,\ldots, X_n\bigr) \rightarrow0
\]
in $P_{G_0}$-probability.
\end{Remark*}

%A simple corollary for consistency.
%
%Let $G_0 \in\Gcal_k(\Theta) \subseteq\Gcal$ for some $k < \infty$,
%and the family of likelihood functions is finitely identifiable.
%Let $\varepsilon> 0$. Suppose that there are sets $\Gcal_n \subset\Gbar(
%and a constant $C>0$ such that
%Moreover, suppose $M_n(\varepsilon)$ is nondecreasing sequence dependent
%on
%$\varepsilon$ so that
%Then,
%$\Pi(G\dvtx  \Wtwo(G_0,G) \geq M_n\varepsilon|X_1,\ldots,X_n)
%}

Before moving to specific examples, we state a simple lemma which
provides estimates of the entropy under the $\Wr$ metric for a number
of classes
of discrete measures of interest. Because $\Wr$ inherits directly
the $\rho$ metric in $\Theta$, the entropy for classes in
$(\Gbar(\Theta), \Wr)$
can typically be bounded in terms of the covering number for subsets of
$(\Theta,\rho)$.
%
%le4 #&#
\begin{lemma}
\label{Lem-entropy}
Suppose that $\Theta$ is bounded. Let $r\geq1$.
\begin{longlist}[(a)]
\item[(a)]
$\log N(2\varepsilon, \Gcal_{k}(\Theta), \Wr) \leq k(\log N(\varepsilon
,\Theta,\rho)
+ \log(e+e\Diam(\Theta)^r/\varepsilon^r))$.\vspace*{1pt}

\item[(b)]
$\log N(2\varepsilon, \bar{\Gcal}(\Theta), \Wr) \leq N(\varepsilon,\Theta,\rho)\log
(e+e\Diam(\Theta)^r/\varepsilon^r)$.

\item[(c)] Let $G_0 = \sum_{i=1}^{k}p^*_i \delta_{\theta^*_i} \in
\Gcal_k(\Theta)$.
Assume that $M = \max_{i=1}^{k}1/p^*_i < \infty$ and
$m = \min_{i,j \leq k} \rho(\theta^*_i,\theta^*_j) > 0$.
Then,
\begin{eqnarray*}
&&
\log N\bigl(\varepsilon/2, \bigl\{G \in\Gcal_k(\Theta)\dvtx
\Wr(G_0,G) \leq 2\varepsilon\bigr\}, \Wr\bigr)
\\
&&\qquad\leq k\Bigl(\sup_{\Theta'} \log N\bigl(\varepsilon/4, \Theta',
\rho\bigr) + \log\bigl(2^{2+3r}k\Diam(\Theta)/m\bigr)\Bigr),
\end{eqnarray*}
where the supremum in the right-hand side is taken over all bounded subsets
\mbox{$\Theta' \subseteq\Theta$} such that $\Diam(\Theta') \leq
4M^{1/r}\varepsilon$.

%is exponentially tight with constant $\eta$, then
%(\log N(\varepsilon,\Theta,\rho) + \log(e+e\Diam(\Theta)/\varepsilon).\]
%is polynomially tight with constant $\eta$, then
%(\log N(\varepsilon,\Theta,\rho) + \log(e+e\Diam(\Theta)/\varepsilon).\]
%}
\end{longlist}
\end{lemma}

%s4 #&#
\section{Examples}
\label{Sec-examples}

In this section we derive posterior contraction rates for
two classes of mixture models, for example, finite mixtures
of multivariate distributions and infinite mixtures based on the
Dirichlet process.

%s4.1 #&#
\subsection{Finite mixtures of multivariate distributions}
Let $\Theta$ be a subset of $\real^d$, $\rho$ be the Euclidean
metric, and $\Pi$ is a prior distribution for discrete measures
in $\Gcal_k(\Theta)$, where $k < \infty$ is known.
Suppose that the ``truth''
$G_0 = \sum_{i=1}^{k}p^*_i \delta_{\theta^*_i} \in\Gcal_k(\Theta)$.
To obtain the convergence rate of the posterior distribution of $G$,
we need the following:

%pa4.1.subsubsection.1 #&#
\subsubsection*{Assumptions A}
(A1) $\Theta$ is compact and the family of likelihood
functions $f(\cdot|\theta)$ is strongly identifiable.

(A2) For some positive constant $C_1$,
$\dK(f_i,f'_j) \leq C_1 \|\theta_i-\theta'_j\|^2$
%and $\int f_i[\log(f_i/f'_j)]^2 \leq C_2 \rho^2(\theta_i,\theta'_j)$
for any \mbox{$\theta_i,\theta'_j \in\Theta$}.
For any $G \in\supp(\Pi)$,
$\int p_{G_0}(\log(p_{G_0}/p_{G}))^2 < C_2 \dK(p_{G_0},p_{G})$
for some constant $C_2>0$.

(A3) Under prior $\Pi$, for small $\delta> 0$,
$c_3 \delta^{k} \leq\Pi(|p_i-p^*_i|\leq\delta, i=1,\ldots, k)
\leq
C_3\delta^{k}$ and
$c_3 \delta^{kd} \leq\Pi(\|\theta_i-\theta^*_i\| \leq\delta,
i=1,\ldots, k) \leq
C_3\delta^{kd}$ for some constants $c_3, C_3 > 0$.

(A4) Under prior $\Pi$, all $p_i$ are bounded away from 0,
and all pairwise distances $\|\theta_i- \theta_j\|$ are bounded away
from 0.

%pa4.1.subsubsection.2 #&#
\begin{Remark*}
Assumptions
(A1) and (A2) hold for the family of Gaussian densities
with mean parameter $\theta$.
Assumption (A3) holds when the prior distribution on the relevant
parameters behaves like
a uniform distribution, up to a multiplicative constant.
\end{Remark*}
%
%th5 #&#
\begin{theorem}
\label{Thm-finite}
Under assumptions \textup{(A1)--(A4)}, the contraction rate in the $L_2$
Wasserstein distance
metric of the posterior distribution of $G$ is $(\log n)^{1/2}n^{-1/4}$.
\end{theorem}
\begin{pf}
Let $G=\sum_{i=1}^{k}p_i \delta_{\theta_i}$.
Combining Lemma~\ref{Lem-gen-ineq} with assumption (A2),
if $\|\theta_i - \theta^*_i\| \leq\varepsilon$
and $|p_i-p^*_i| \leq\varepsilon^2/(k\Diam(\Theta)^2)$ for
$i=1,\ldots, k$, then
$\dK(p_{G_0},\break p_{G}) \leq d_{\rho_K}(G_0,G)
\leq C_1 \sum_{1\leq i,j \leq k} q_{ij} \|\theta^*_i-\theta_j\|^2$,
for any $\vecq\in\Qcal$. Thus,
$\dK(p_{G_0},\break p_{G}) \leq C_1 \Wtwo^2(G_0,G) \leq
C_1 \sum_{i=1}^{k}(p^*_i \wedge p_i) \|\theta^*_i-\theta_i\|^2 +\vadjust{\goodbreak}
C_1 \sum_{i=1}^{k}|p_i-\break p^*_i|\Diam(\Theta)^2 \leq2C_1 \varepsilon^2$.
Hence, under prior $\Pi$,
\begin{eqnarray*}
&&\Pi\bigl(G\dvtx  \dK(p_{G_0},p_{G}) \leq\varepsilon^2
\bigr)
\\
&&\qquad\geq\Pi\bigl(G\dvtx  \bigl\|\theta_i-\theta^*_i\bigr\| \leq\varepsilon,
\bigl|p_i-p^*_i\bigr|\leq \varepsilon^2/\bigl(k\Diam(
\Theta)^2\bigr), i=1,\ldots, k\bigr). %\geq c_3^2 \varepsilon^{kd}(\varepsilon^2/(k\Diam(\Theta)^2))^{k}
\end{eqnarray*}
%
%Similarly, due to (A2), $\int p_{G_0}[\log(p_{G_0}/p_G)]^2 \leq C_2d_{
In view of assumptions (A2) and (A3), we have $\Pi(B_K(\varepsilon))
\geqa\varepsilon^{k(d+2)}$.
Conversely, for sufficiently small $\varepsilon$, if $\Wtwo(G_0,G) \leq
\varepsilon$,
then by reordering the index of the atoms, we must have $\|\theta_i-\theta^*_i\|
= O(\varepsilon)$ and $|p_i-p^*_i|= O(\varepsilon^2)$ for all
$i=1,\ldots, k$ [see
the argument in the proof of Lemma~\ref{Lem-entropy}(c)]. This
entails that under the prior $\Pi$,
\begin{eqnarray*}
\Pi\bigl(G\dvtx  \Wtwo^2(G_0,G) \leq\varepsilon^2
\bigr) & \leq& \Pi\bigl(G\dvtx  \bigl\|\theta_i-\theta^*_i\bigr\| \leq
O(\varepsilon), \bigl|p_i-p^*_i\bigr|\leq O\bigl(
\varepsilon^2\bigr), \forall i\bigr)
\\
&\lesssim& \varepsilon^{k(d+2)}.
\end{eqnarray*}
Let $\varepsilon_n$ be a sufficiently large multiple of $(\log n/n)^{1/2}$.
We proceed by verifying conditions of Theorem
\ref{Thm-convergence-2}.
Let $\Gcal_n:= \Gcal_k(\Theta)$. Then $\Pi(\Gbar(\Theta
)\setminus\Gcal_n) = 0$,
so equation (\ref{Eqn-support-1}) trivially holds.

Next, we provide upper bounds for $D(\varepsilon/2, S, \Wtwo)$, where
$S = \{G\in\Gcal_n\dvtx \break \Wtwo(G_0, G) \leq2\varepsilon\}$,
and $M(\Gcal_n,G_1,\varepsilon)$ so that (\ref{Eqn-entropy-2}) is
satisfied. Indeed, for any $\varepsilon> 0$, $\log D(\varepsilon/2, S,
\Wtwo) \leq\log N(\varepsilon/4, S, \Wtwo)$. By Lemma
\ref{Lem-entropy}(c) and assumption (A4), $N(\varepsilon/4, S, W_1)$ is
bounded in terms of $\sup_{\Theta'} \log N(\varepsilon/8, \Theta',
\rho)$, which is bounded above by a constant when $\Theta'$'s are
subsets of $\Theta$ whose diameter is bounded by a multiple of
$\varepsilon$. Turning to $M(\Gcal_n, G_1,\varepsilon)$, due to strong
identifiability and assumption (A2), $\Psi_{\Gcal_n}(\varepsilon) \geq
c\varepsilon^4$ for some constant $c>0$. By Lemma~\ref{Lem-entropy}(a),
for some constant $c_1 > 0$, $\log M(\Gcal_n,G_1,\varepsilon) \leq\log
N(c_1\varepsilon^2, \Gcal_k(\Theta) \cap B_W(G_1,\varepsilon /2),
\Wtwo) \leq k(\log N(c_1\varepsilon^2/2, \Theta, \rho) + \log(e+4e\Diam
(\Theta)^2/ c_1^2\varepsilon^4)) \leq \break n\varepsilon_n^2/2$. Thus,
equation (\ref{Eqn-entropy-2}) holds.

By Proposition~\ref{Thm-1}(b)
and assumption (A4), we have
\[
\Psi_{\Gcal_n}(j\varepsilon_n) = \inf_{\Wtwo(G_0,G) \geq j\varepsilon/2}
\dHel^2(p_{G_0},p_{G}) \geq c(j
\varepsilon_n)^{4}
\]
for some constant $c> 0$.
To ensure condition (\ref{Eqn-rate-cond}),
note that (constants $c$ change after each bounding step)
\begin{eqnarray*}
\exp\bigl(2n\varepsilon_n^2\bigr)\sum
_{j\geq M_n} \exp\bigl[-n\Psi_{\Gcal
_n}(j\varepsilon_n)/16
\bigr] &\lesssim&\exp\bigl(2n\varepsilon_n^2\bigr)\sum
_{j\geq M_n} \exp\bigl[-nc(j\varepsilon_n)^{4}
\bigr]
\\
&\leqa& \exp\bigl(2n\varepsilon_n^2 -ncM_n^4
\varepsilon_n^{4}\bigr).
\end{eqnarray*}
This upper bound goes to zero if $ncM_n^4\varepsilon_n^4 \geq4n\varepsilon_n^2$,
which is satisfied by taking $M_n$ to be a large multiple of $\varepsilon_n^{-1/2}$.
Thus, we need $M_n\varepsilon_n \asymp\varepsilon_n^{1/2} \asymp(\log n)^{1/4}
n^{-1/4}$.

Under the assumptions
specified above,
\[
\Pi\bigl(G\dvtx  j\varepsilon_n < \Wtwo(G,G_0) \leq2j
\varepsilon_n\bigr)/ \Pi\bigl(B_K(\varepsilon_n)
\bigr) = O(1).
\]
On the other hand, for $j \geq M_n$,
we have $\exp[n\Psi_{\Gcal_n}(j\varepsilon_n)/16] \geq
\exp[nc(j\varepsilon_n)^4/16]$ which is bounded below by an arbitrarily
large constant
by choosing $M_n$ to be a large multiple of $\varepsilon_n^{-1/2}$,
thereby ensuring (\ref{Eqn-support-2}).\vadjust{\goodbreak}

Thus, by Theorem~\ref{Thm-convergence-2}, rate of contraction for the posterior
distribution of $G$ under the $\Wtwo$ distance metric is $(\log
n)^{1/4}n^{-1/4}$,
which is up to a logarithmic factor the minimax optimal rate
$n^{-1/4}$ as proved for the univariate finite mixtures
by~\cite{Chen-95}.\vspace*{-3pt}
\end{pf}

%s4.2 #&#
\subsection{Dirichlet process mixtures}

Given the ``true'' discrete mixing measure,
$G_0 = \sum_{i=1}^{k}p^*_i \delta_{\theta^*_i}
\in\Gcal_{k}(\Theta)$, where $\Theta$ is a metric space
but $k \leq\infty$ is unknown.
To estimate $G_0$, the prior distribution $\Pi$ on discrete measure
$G\in\Gbar(\Theta)$ is taken to be a Dirichlet process $\operatorname{DP}(\nu,P_0)$
that centers at $P_0$ with concentration parameter $\nu> 0$ \cite
{Ferguson-73}.
Here, parameter $P_0$ is a probability measure on $\Theta$.
For any $r \geq1$,
the following lemma provides a lower bound of small ball probabilities
of metric
space $(\Gbar(\Theta), \Wr)$ in terms of
small ball $P_0$-probabilities of metric space
$(\Theta,\rho)$.\vspace*{-3pt}
%
%le5 #&#
\begin{lemma}
\label{Lem-DPdense}
Let $G\sim\operatorname{DP}(\nu,P_0)$, where $P_0$ is a nonatomic base
probability measure
on a compact set $\Theta$. For a small $\varepsilon> 0$, let
$D = D(\varepsilon, \Theta, \rho)$ denote the $\varepsilon$-packing
number of $\Theta$ under the
$\rho$ metric.
Then, under the Dirichlet process distribution,
\begin{eqnarray*}
\Pi\bigl(G\dvtx  \Wr^r(G_0,G) \leq\bigl(2^\m+1
\bigr)\varepsilon^\m\bigr)
\geq\frac{\Gamma(\nu)\nu^D}{(2D)^{D-1}} \biggl(\frac{\varepsilon}{\Diam(\Theta)} \biggr)^{\m(D-1)}
\sup_{S}\prod_{i=1}^{D}P_0(S_i).
\end{eqnarray*}
%
%In addition, if $D \leqa[\Diam(\Theta)/\varepsilon]^\gamma$ for some $
Here, $S:=(S_1,\ldots,S_D)$ denotes the $D$ disjoint $\varepsilon
/2$-balls that
form a maximal $\varepsilon$-packing of $\Theta$. The supremum is taken
over all such packings. $\Gamma(\cdot)$ is the gamma
function.\vspace*{-2pt}
\end{lemma}
\begin{pf}
Since every point in $\Theta$ is of distance at most $\varepsilon$ to
one of the
centers of $S_1,\ldots,S_D$, there is a $D$-partition $(S'_1,\ldots,S'_D)$ of
$\Theta$, such that $S_i \subseteq S'_i$, and $\Diam(S'_i) \leq
2\varepsilon$
for each $i=1,\ldots, D$. Let $m_i = G(S'_i)$, $\mu_i = P_0(S'_i)$,
and $\hatp_i = G_0(S'_i)$. From the definition of Dirichlet processes,
$\mathbf{m} = (m_1,\ldots, m_D) \sim\operatorname{Dir}(\nu\mu_1,\ldots, \nu\mu_{D})$.
To obtain an upper bound for $d_{\rho^\m}(G_0,G)$, consider a coupling
between $G_0$ and $G$, by associating $m_i \wedge\hat{p}_i$ probability
mass of supporting atoms for $G_0$ contained in subset $S'_i$ with the
same probability
mass of supporting atoms for $G$ contained in the same subset, for each
$i=1,\ldots,D$.
The remaining mass (of probability $\|\mathbf{m}-\hat{\mathbf{p}}\|
$) for both measures
are coupled in an arbitrary way. The expectation under this coupling
of the $\rho^\m$ distance provides one such upper bound, that is,
\[
d_{\rho^\m}(G_0,G) \leq(2\varepsilon)^\m+ \|
\mathbf{m} - \hat {\mathbf{p}}\|_1 \bigl[\Diam(\Theta)
\bigr]^\m.
\]

Due to the nonatomicity of $P_0$, for $\varepsilon$ sufficiently small,
$\nu\mu_i \leq1$ for all $i=1,\ldots, D$. Let $\delta= \varepsilon
/\Diam(\Theta)$.
Then, under $\Pi$,
\begin{eqnarray*}
&&
\operatorname{Pr}\bigl(d_{\rho^\m}(G_0,G) \leq
\bigl(2^\m+1\bigr)\varepsilon^\m\bigr)\\
&&\qquad\geq \operatorname{Pr}\bigl(
\|\mathbf{m}-\hat{\mathbf{p}}\|_1 \leq\delta^\m\bigr)
\\
&&\qquad \geq\operatorname{Pr}\bigl(|m_i-\hatp_i| \leq
\delta^\m/2D, i=1,\ldots,D-1\bigr)
\\
&&\qquad = \frac{\Gamma(\nu)}{\prod_{i=1}^{D}\Gamma(\nu\mu_i)} \int_{\Delta_{D-1} \cap|m_i-\hatp_i|\leq\delta^\m/2D} \prod
_{i=1}^{D-1}m_i^{\nu\mu_i-1}\Biggl(1-
\sum_{i=1}^{D-1}m_i
\Biggr)^{\nu\mu_D-1} \,d m_i
\\
&&\qquad \geq\frac{\Gamma(\nu)}{\prod_{i=1}^{D}\Gamma(\nu\mu_i)} \prod_{i=1}^{D-1}
\int_{\max(\hatp_i-\delta^\m/2D,0)}^{\min
(\hatp_i+\delta^\m/2D,1)} m_i^{\nu\mu_i-1} \,d
m_i
\\
&&\qquad \geq\Gamma(\nu) \bigl(\delta^\m/2D\bigr)^{D-1}\prod
_{i=1}^{D}(\nu\mu_i).
\end{eqnarray*}
The second inequality in the previous display is due to the fact that
$\|\mathbf{m}-\hat{\mathbf{p}}\|_1 \leq2 \sum_{i=1}^{D-1} |m_i -
\hat{p}_i|$.
The third inequality is due to $(1-\sum_{i=1}^{D-1}m_i)^{\nu\mu_D-1}
= m_{D}^{\nu\mu_D -1} \geq1$,
since $\nu\mu_D \leq1$ and $0 < m_{D} <1$ almost surely.
The last inequality is due to the fact that
$\Gamma(\alpha) \leq1/\alpha$ for $0< \alpha\leq1$. This gives
the desired claim.
\end{pf}

%pa4.2.subsubsection.1 #&#
\subsubsection*{Assumptions B}

(B1) The nonatomic base measure $P_0$ places full support on
a bounded set
$\Theta\subset\real^d$. Moreover,
$P_0$ has a Lebesgue density that is bounded away from zero.
%The family of the likelihood densities $f(\cdot|\theta)$ is finitely
%identifiable.

(B2) For some constants $C_1,m_1 >0$,
$\dK(f_i,f'_j) \leq C_1 \rho^{m_1}(\theta_i,\theta'_j)$
%$\int f_i [\log(f_i/f'_j)]^2 \leq C_1 \rho^{m_2}(\theta_i,\theta'_j)$
for any \mbox{$\theta_i,\theta'_j \in\Theta$}.

For any $G \in\supp(\Pi)$,
$\int p_{G_0}(\log(p_{G_0}/p_{G}))^2 \leq C_2 \dK(p_{G_0},p_{G})^{m_2}$
for some constants $C_2,m_2>0$.
%
%$P_0$ places sufficient probability mass on
%all small balls that pack $\Theta$.
%Specifically, there is a universal constant $c_3>0$ such
%that the probability of the
%$D$-partition $(S_1,\ldots, S_D)$ specified in Lemma
%
%so that the packing number
%$D(\varepsilon,\Theta,\rho) \asymp[\Diam(\Theta)/\varepsilon]^d$.
%
%$D(\varepsilon,\Theta,\rho) \asymp\exp[(1/\varepsilon)^{1/\gamma}]$.
%
%th6 #&#
\begin{theorem}
Given assumptions \textup{(B1)} and \textup{(B2)}
and the smoothness conditions for the likelihood family
as specified in Theorem~\ref{Thm-convolution},
there is a sequence $\beta_n \searrow0$
such that $\Pi(\Wtwo(G_0,G) \geq\beta_n |X_1,\ldots,X_n)
\rightarrow0$
in $P_{G_0}$ probability. Specifically:
\begin{longlist}[(2)]
\item[(1)] For ordinary smooth likelihood functions, take
\[
\beta_n \asymp(\log n/n)^{{2}/({(d+2)(4+(2\beta+1)d')})}
\]
for any constant $d' > d$.

\item[(2)] For supersmooth likelihood functions, take
$\beta_n \asymp(\log n)^{-1/\beta}$.
\end{longlist}
\end{theorem}
\begin{pf}
The proof consists of two main steps.
First, we shall prove that under assumptions (B1)--(B2), conditions
specified by
(\ref{Eqn-entropy-2s}), (\ref{Eqn-support-1s}) and (\ref
{Eqn-support-2s}) in
Theorem~\ref{Thm-convergence-1} are satisfied by taking $\Gcal_n =
\Gbar(\Theta)$,
which is a convex set, and $\varepsilon_n$ to be a large multiple of
$(\log n/n)^{1/(d+2)}$.
%If Assumption (B4) is replaced by (B'4), then the same claim holds
%by
%taking $\varepsilon_n$ to be a large multiple of
%$1/(\log n)^{\gamma}$.
The second step involves constructing a sequence of $M_n$ and
$\beta_n = M_n\varepsilon_n$ for which Theorem~\ref{Thm-convergence-1} can
be applied.

\textit{Step} 1: By Lemma~\ref{Lem-gen-ineq} and (B2),
$\dK(p_{G_0},p_{G}) \leq d_{\rho_K}(G_0,G) \leq C_1 d_{\rho^{m_1}}(G_0,G)$.
Also, $\int p_{G_0}[\log(p_{G_0}/p_{G})]^2 \leqa C_2 d_{\rho
^{m_1}}(G_0,G)^{m_2}$.
Assume without loss of generality that $m_1 \leq m_2$ (the other direction
is handled similarly).\vadjust{\goodbreak}
We obtain that
$\Pi(G\in B_K(\varepsilon_n)) \geq
\Pi(G\dvtx  d_{\rho^{m_1}}(G_0,G) \leq C_3 \varepsilon_n^{2 \vee2/m_2})$
for some constant $C_3 > 0$.

From (B1), there is a universal constant $c_3>0$ such that for any
$\varepsilon$
and any $D$-partition $(S_1,\ldots, S_D)$ specified in Lemma
\ref{Lem-DPdense}, there holds
\[
\log\prod_{i=1}^{D}P_0(S_i)
\geq c_3D\log(1/D).
\]
Moreover, the packing number satisfies $D \asymp[\Diam(\Theta
)/\varepsilon_n]^d$.
Combining these facts with Lemma~\ref{Lem-DPdense},
we have
$\log\Pi(G\in B_K(\varepsilon_n)) \geqa(D-1)\log(\varepsilon_n/\Diam
(\Theta))
+ (2D-1)\log(1/D) + D\log\nu$, where the approximation constant is
dependent on $m_1,m_2$.
It is simple to check that
condition (\ref{Eqn-support-2s}) holds,
$\log\Pi(G\in B_K(\varepsilon_n)) \geq-Cn\varepsilon_n^2$, by the given
rate of $\varepsilon_n$,
for any constant $C>0$.\vspace*{2pt}

Since $\Gcal_n = \Gbar(\Theta)$,\vspace*{1pt} (\ref{Eqn-support-1s}) trivially holds.
Turning to condition (\ref{Eqn-entropy-2s}), by Lem\-ma~\ref
{Lem-entropy}(b), we
have $\log N(2\varepsilon_n,\Gbar(\Theta),\Wtwo)
\leq N(\varepsilon_n, \Theta, \rho)
\log(e + e \Diam(\Theta)^2/\varepsilon_n^2)
\leq(\Diam(\Theta)/\varepsilon_n)^{d}
\log(e + e \Diam(\Theta)^2/\varepsilon_n^2) \leq n \varepsilon_n^2$ by
the specified
rate of $\varepsilon_n$.

%If (B5) is replaced by (B5'), again, similar calculations show that
%all conditions of Theorem~\ref{Thm-convergence-1} holds
%if $\varepsilon_n$ is a large multiple of $1/(\log n)^{\gamma}$.

\textit{Step} 2: For any $\Gcal\subseteq\Gbar(\Theta)$, let
$R_{\Gcal}(r)$ be the
inverse of the Hellinger information function of the $W_2$ metric. Specifically,
for any $t \geq0$,
\[
R_{\Gcal}(t) = \inf\bigl\{r \geq0 | \Psi_\Gcal(r) \geq t\bigr
\}.
\]
Note that $R_{\Gcal}(0) = 0$. $R_{\Gcal}(\cdot)$ is nondecreasing because
$\Psi_{\Gcal}(\cdot)$ is.
%Moreover, $R_{\Gbar(\Theta)}(t) \searrow0$ as $t \rightarrow0$.
%Indeed,
%if this is not true, then there is a sequence of $t_m \rightarrow0$
%and
%$\delta> 0$ such that $R_{\Gbar(\Theta)}(t_m) > \delta$ for all $m
%This implies that $\Psi_{\Gbar(\Theta)}(r) = 0$ for any $r < \delta$,
%which
%leads to contradiction by Proposition~\ref{Thm-1}(a), due to the
%compactness of
%$\Gbar(\Theta)$ and finite identifiability from (B1).
%}

Let $(\varepsilon_n)_{n\geq1}$ be the sequence determined in the previous
step of the proof. Let $M_n = R_{\Gbar(\Theta)}(8\varepsilon_n^2(C+4))/\varepsilon_n$,
and $\beta_n = M_n \varepsilon_n = R_{\Gbar(\Theta)}(8\varepsilon_n^2(C+4))$.
Condition (\ref{Eqn-rate-conds-2s}) holds by definition of $R_{\Gbar
(\Theta)}$,
that is, $\Psi_{\Gcal(\Theta)}(M_n \varepsilon_n) \geq8\varepsilon_n^2(C+4)$.
To verify (\ref{Eqn-hellinger-sum}), note that the running
sum with respect to $j$ cannot have more than $\Diam(\Theta)/\varepsilon_n$ terms,
and, due to the monotonicity of $\Psi_{\Gcal}$, we have
\begin{eqnarray*}
&&
\exp\bigl(2n\varepsilon_n^2\bigr)\sum
_{j\geq M_n} \exp\bigl[-n\Psi_{\Gcal
_n}(j\varepsilon_n)/8
\bigr]
\\
&&\qquad\leq\Diam(\Theta)/\varepsilon_n \exp\bigl(2n\varepsilon_n^2
- n\Psi_{\Gcal
_n}(M_n\varepsilon_n)/8\bigr)
\rightarrow0.
\end{eqnarray*}
Hence, Theorem~\ref{Thm-convergence-1} can be applied to
conclude that
\[
\Pi\bigl(W_2(G_0,G) \geq\beta_n |
X_1,\ldots,X_n\bigr) \rightarrow0
\]
in $P_{G_0}$-probability.
Under the ordinary smoothness condition
(as specified in Theorem~\ref{Thm-convolution}),
$R_{\Gbar(\Theta)}(t) = t^{{1}/({4+(2\beta+1)d + \delta})}$, where
$\delta$ is an arbitrarily positive constant. So,
\[
\beta_n \asymp\varepsilon_n^{{2}/({4+(2\beta+1)d +\delta})} = (\log
n/n)^{{2}/({(d+2)(4+(2\beta+1)d +\delta)})}.
\]
On the other hand, under the supersmoothness condition,
$R_{\Gbar(\Theta)}(t) = (1/\break\log(1/ t))^{1/\beta}$.
So, $\beta_n \asymp(\log(1/\varepsilon_n))^{-1/\beta}
\asymp(\log n)^{-1/\beta}$.
\end{pf}

\section{Proofs}
\label{Sec-proofs}

%s5.1 #&#
\subsection{Proofs of Wasserstein identifiability results}

%Let $G,G' \in\Gbar(\Theta)$ such that $\drPhi(G,G') < \infty$ for
%some convex function $\phi$. \\
%Then, $\dPhi(p_G,p_{G'}) \leq\drPhi(G,G')$.
%%In particular, $d_{h}^2(p_G,p_{G'}) \leq d_{\rho h^2}(G,G')$,
%%$d_{V}(p_G,p_{G'}) \leq d_{\rho V} (G,G')$, and
%%$\dK(p_G,p_{G'}) \leq d_{\rho_K}(G,G')$.
%}

%Suppose that $\Theta$ is compact, the family $\{f(\cdot|\theta),\theta
%is strongly identifiable, and for all $x \in\Xcal$, the Hessian
%matrix $D^{2}f(x|\theta)$ satisfies a uniform Lipschitz condition
%|\gamma^T(D^{2}f(x|\theta_1) - D^{2}f(x|\theta_2))\gamma|
%for all $x, \theta_1,\theta_2$ and some fixed $C$ and $\delta> 0$.
%Then,
%for fixed $G_0 \in\Gcal_k(\Theta)$, where $k<\infty$:
%d_\rho(G_0,G') \leq\varepsilon\biggr\} > 0.
%The assertion also holds with $\psi$ being replaced by $\psi_1$.

\mbox{}

\begin{pf*}{Proof of Theorem~\ref{Thm-bound-supnorm}}
Suppose that equation (\ref{Eqn-bound-supnorm}) is not true, then
there will be
sequences of $G_{n}$ and $G'_{n}$ tending to $G_0$ in the $\Wtwo$ metric,
and that \mbox{$\psi(G_n,G'_n) \rightarrow0$}. We write
$G_n = \sum_{i=1}^{\infty}\p_{n,i} \delta_{\theta_{n,i}}$, where
$p_{n,i}=0$
for indices $i$ greater than $k_n$, the number of atoms of $G_n$.
Similar notation is
applied to $G'_n$.
Since both $G_n$ and $G'_n$ have a finite number of atoms,
there is $\qvec^{(n)} \in\Qcal(\vecp_n, \mathbf{p}'_n)$ so that
$\Wtwo^2(G_n,G'_n) = \sum_{ij}q_{ij}^{(n)}\|\theta_{n,i}-\theta'_{n,j}\|^2$.
%Note that $d_\rho^2(G_n,G'_n) \leq d_{\rho^2}(G_n,G'_n) = O(d_
%while the latter inequality is due to the boundedness of $\Theta$.

Let $\Ocal_n = \{(i,j)\dvtx  \|\theta_{n,i}-\theta'_{n,j}\|
\leq\Wtwo(G_n,G'_n)\}$. This set exists because there are
pairs of atoms $\theta_{n,i}, \theta'_{n,j}$ such that
$\|\theta_{n,i}-\theta'_{n,j}\|$ is bounded away from zero in the limit.
%Then, $\sum_{(i,j)\not\in\Ocal_n}
%q_{ij}^{(n)} \leq d_{\rho^2}(G_n,G'_n)/d_{\rho^2}^{1-\delta'}(G_n,G'_n)
%
Since $\qvec^{(n)} \in\Qcal(\vecp_n, \mathbf{p}'_n)$, we can express
\begin{eqnarray*}
\psi\bigl(G_n,G'_n\bigr) & = &
\sup_{x} \Biggl| \sum_{i=1}^{k_n}
\p_{n,i} f(x|\theta_{n,i}) - \sum
_{j=1}^{k'_n} \p'_{n,j} f
\bigl(x|\theta'_{n,j}\bigr) \Biggr| \Big/\Wtwo^2
\bigl(G_n,G'_n\bigr)
\\
& = & \sup_{x} \biggl| \sum_{ij}
q_{ij}^{(n)}\bigl(f(x|\theta_{n,i})-f\bigl(x|
\theta'_{n,j}\bigr)\bigr) \biggr| \Big/\Wtwo^2
\bigl(G_n,G'_n\bigr)
\end{eqnarray*}
and, by Taylor's expansion,
\begin{eqnarray*}
&&
\psi\bigl(G_n,G'_n\bigr) \\
&&\qquad =
\sup_{x} \biggl| \sum_{(i,j) \notin\Ocal
_n}
q_{ij}^{(n)}\bigl(f\bigl(x|\theta'_{n,j}
\bigr)-f(x|\theta_{n,i})\bigr)
\\
&&\qquad\quad\hspace*{17pt}{} +\sum_{(i,j) \in\Ocal_n} q_{ij}^{(n)}
\bigl(\theta'_{n,j}-\theta_{n,i}
\bigr)^T Df(x|\theta_{n,i})
\\
&&\qquad\quad\hspace*{17pt}{} +\sum_{(i,j) \in\Ocal_n} q_{ij}^{(n)}
\bigl(\theta'_{n,j}-\theta_{n,i}
\bigr)^T D^{2}f(x|\theta_{n,i}) \bigl(
\theta'_{n,j}-\theta_{n,i}\bigr) + R_n(x) \biggr|
\\
&&\qquad\quad{} \Big/\Wtwo^2\bigl(G_n,G'_n
\bigr)
\\
&&\qquad =: \sup_{x}\bigl|A_n(x) + B_n(x) +
C_n(x) + R_n(x)\bigr|/D_n,
\end{eqnarray*}
where
\[
R_n(x) = O \biggl(\sum_{(i,j)\in\Ocal_n}
q_{ij}^{(n)}\bigl\|\theta_{n,i}-\theta'_{n,j}
\bigr\|^{2+\delta} \biggr) = o \biggl(\sum_{(i,j)\in\Ocal_n}
q_{ij}^{(n)}\bigl\|\theta_{n,i}-\theta'_{n,j}
\bigr\|^2 \biggr)
\]
due to (\ref{Eqn-Lipschitz}) and the definition of $\Ocal_n$.
So $R_n(x)/\Wtwo^2(G_n,G'_n) \rightarrow0$.

The quantities $A_n(x), B_n(x)$
and $C_n(x)$ are linear combinations of elements of
$f(x|\theta)$, $Df(x|\theta)$
and $D^{2}f(x|\theta)$ for different $\theta$'s, respectively.
Since $\Theta$ is compact,
subsequences\vadjust{\goodbreak} of $G_n$ and $G'_n$ can be chosen so that each
of their support points converges to a fixed atom $\theta^*_l$, for
$l=1,\ldots,k^* \leq k$.
After being rescaled, the limits of $A_n(x)/D_n, B_n(x)/D_n$ and $C_n(x)/D_n$
are still linear combinations with constant coefficients not depending
on $x$.

We shall now argue that not all such coefficients vanish to zero.
Suppose this is not the case. It follows that for
the coefficients of $C_n(x)/D_n$ we have
\[
\sum_{(i,j)\in\Ocal_n} q_{ij}^{(n)}\bigl\|
\theta'_{n,j}-\theta_{n,i}\bigr\|^2/
\Wtwo^2\bigl(G_n,G'_n\bigr)
\rightarrow0.
\]
This implies that $\sum_{(i,j)\notin\Ocal_n} q_{ij}^{(n)}\|\theta'_{n,j}-\theta_{n,i}\|^2/
\Wtwo^2(G_n,G'_n) \rightarrow1$. Since $\Theta$ is bounded,
there exists a pair $(i,j) \notin\Ocal_n$ such that
$q_{ij}^{n}/\Wtwo^2(G_n,G'_n)$ does not vanish to zero. But then,
one of the coefficients of $A_n(x)/D_n$ does not vanish to zero,
which contradicts the hypothesis.

Next, we observe that some of the coefficients of $A_n(x)/D_n,
B_n(x)/D_n$ and
$C_n(x)/D_n$ may tend to infinity.
For each $n$, let $d_n$ be the inverse of the maximum coefficient of
$A_n(x)/D_n$, $B_n(x)/D_n$, and $C_n(x)/D_n$. From the conclusion in
the preceding paragraph, $|d_n|$ is uniformly bounded from above
by a constant for all $n$.
Moreover, $d_n A_n(x)/D_n$ converges to $\sum_{j=1}^{k^*}\alpha_{j}f(x|\theta^*_{j})$ and $d_n B_n(x)/D_n$ converges to $\sum_{j=1}^{k^*}\beta_{j}^TDf(x|\theta^*_{j})$,
and $d_n C_n(x)/D_n$ converges to $\sum_{j=1}^{k^*}\gamma_{j}D^2f(x|\theta^*_{j})
\gamma_{j}$, for some finite $\alpha_{j}, \beta_{j}$ and $\gamma_{j}$,
not all of them vanishing (in fact, at least one of them is 1).
We have
%
%e23 #&#
\begin{eqnarray}\label{Eqn-conv-pointwise}
&&
d_n \bigl|p_{G_n}(x)-p_{G'_n}(x)\bigr|/\Wtwo^2
\bigl(G_n,G'_n\bigr)
\nonumber\\[-8pt]\\[-8pt]
&&\qquad\rightarrow\Biggl|\sum_{j=1}^{k^*}\alpha_{j} f
\bigl(x|\theta^*_{j}\bigr) + \beta_j^T Df
\bigl(x|\theta^*_j\bigr) + \gamma_j^T
D^{2}f\bigl(x|\theta^*_j\bigr) \gamma_j \Biggr|
\nonumber
\end{eqnarray}
for all $x$. This entails that the right-hand side of the preceding display
must be 0 for almost all $x$. By strong identifiability, all
coefficients must be 0, which leads to contradiction.

With respect to $\psi_1(G,G')$, suppose that the claim is not true,
which implies the existence of a subsequence $G_n,G'_n$ such that
that $\psi_1(G_n,G'_n) \rightarrow0$. Going through the same argument
as above,
we have $\alpha_{j}, \beta_{j}, \gamma_{j}$, not all of which are
zero, such
that equation (\ref{Eqn-conv-pointwise}) holds. An application of
Fatou's lemma
yields $\int|\sum_{j=1}^{k^*}\alpha_{j} f(x|\theta_{j}) +
\beta_j^T Df(x|\theta_j) + \gamma_j^T D^{2}f(x|\theta_j) \gamma_j
| \,d\mu= 0$.
Thus, the integrand must be 0 for almost all $x$, leading to
contradiction.~%
\end{pf*}
\begin{pf*}{Proof of Theorem~\ref{Thm-convolution}}
To obtain an upper bound of $\Wtwo^2(G,G') =
d_{\rho^2}(G, G')$ in terms of $\dV(p_G,p_{G'})$
under the condition that $\dV(p_G,p_{G'}) \rightarrow0$,
our strategy is approximate $G$ and $G'$ by convolving these
with some mollifier $K_\delta$. By the triangular inequality, $\Wtwo(G,G')
\leq\Wtwo(G,G*K_\delta) + \Wtwo(G',G'*K_\delta) +
\Wtwo(G*K_\delta,G'*K_\delta)$. The first
two terms are simple to bound, while the last term can be handled by expressing
$G*K_\delta$ as the convolution of the mixture density $p_G$ with another
function.
%This trick was widely exploited in kernel density estimation
%method for deconvolution problems (e.g.,~\cite{Zhang-90,Fan-91}).
We also need the following elementary lemma (whose proof is given
Section~\ref{sec-aux}).
%
%le6 #&#
\begin{lemma}
\label{Lem-ineq-moment}
Assume that $p$ and $p'$ are two probability density functions on
$\real^d$
with bounded $\s$-moments.
\begin{longlist}[(a)]
\item[(a)] For $\kk$ such that $0<\kk<\s$,
\[
\int\bigl|p(x)-p'(x)\bigr|\|x\|^\kk dx \leq 2
\bigl\|p-p'\bigr\|_{L_1}^{(\s-\kk)/\s}\bigl(\E_p\|X
\|^\s+ \E_{p'}\|X\|^\s \bigr)^{\kk/\s}.
\]
\item[(b)] Let $V_d= \pi^{d/2}\Gamma(d/2+1)$ denote the volume of
the $d$-dimensional
unit sphere. Then,
\[
\bigl\|p-p'\bigr\|_{L_1} \leq2 V_d^{s/(d+2s)}
\bigl(\E_p \|X\|^s+ \E_{p'}\|X
\|^{s}\bigr)^{{d}/({d+2s})} \bigl\|p-p'\bigr\|_{L_2}^{{2s}/({d+2s})}.
\]
\end{longlist}
\end{lemma}

Take any $s>0$, and let $K\dvtx \real^d \rightarrow(0,\infty)$ be a symmetric
density function on $\real^d$ whose Fourier transform $\tilde{K}$ is
a continuous
function whose support is bounded in $[-1,1]^d$. Moreover, $K$ has
bounded moments
up to order $s$. Consider mollifiers
$K_\delta(x) = \frac{1}{\delta^d} K(x/\delta)$ for $\delta> 0$.
Let $\tilde{K}_\delta$ and $\tilde{f}$ be the Fourier
transforms for $K_\delta$ and~$f$, respectively. Define $g_\delta$ to
be the
inverse Fourier transform of $\tilde{K}_\delta/ \tilde{f}$:
\[
g_\delta(x) = \frac{1}{(2\pi)^d}\int_{\real^d}
e^{i\langle\omega, x \rangle}\frac{\tildeK_\delta(\omega)}{\tildef(\omega)} \,d\omega.
\]

Note that function $\tildeK_\delta(\omega)/\tildef(\omega)$ has
bounded support.
So, $g_\delta\in L_1(\real)$, and $\tildeg_\delta:= \tildeK_\delta
(\omega)/\tildef(\omega)$
is the Fourier transform of $g_\delta$. By the convolution theorem,
$f*g_\delta= K_\delta$. As a result,
\[
G*K_\delta= G*f*g_\delta= p_G *
g_\delta.
\]

From the definition of $K_\delta$, the second moment under $K_\delta$ is
$O(\delta^2)$. Consider a coupling $G$ and
$G*K_\delta$ under which
we have a pair of random variables $(\theta, \theta+\varepsilon)$
where $\varepsilon$ is independent of $\theta$, the marginal
distributions of $\theta$ and $\varepsilon$ are $G$
$K_\delta$, respectively. Under this coupling,
$\E\|(\theta+\varepsilon) - \theta\|^2 = O(\delta^2)$, which
entails that $\Wtwo^2(G,G*K_\delta) = O(\delta^2)$.

By the triangular
inequality, $\Wtwo(G,G') \leq\Wtwo(G*K_\delta,G'*K_\delta) +
O(\delta)$, so for some constant $C(K)>0$ dependent only on kernel $K$,
%
%e24 #&#
\begin{equation}
\label{Eqn-perturb} \Wtwo^2\bigl(G,G'\bigr) \leq2
\Wtwo^2\bigl(G*K_\delta, G'*K_\delta
\bigr) + C(K)\delta^2.
\end{equation}

Theorem 6.15 of~\cite{Villani-09} provides an upper bound for
the Wasserstein distance: for any two probability measures
$\mu$ and $\nu$,
$\Wtwo^2(\mu,\nu) \leq2\int\|x\|^2 \,d|\mu-\nu|(x)$,
where $|\mu-\nu|$ is the total variation of measure $|\mu- \nu|$. Thus,
%
%e25 #&#
\begin{equation}
\label{Eqn-bound-wasser-variational} \Wtwo^2\bigl(G*K_\delta,
G'*K_\delta\bigr) \leq 2 \int\|x\|^2
\bigl|G*K_\delta(x)- G'*K_\delta(x)\bigr| \,dx.
\end{equation}

We note that since density function $K$ has a bounded $s$th moment,
\begin{eqnarray*}
\int\|x\|^s G*K_\delta(dx) & \leq& 2^{s} \biggl[
\int\|\theta\|^s \,dG(\theta) + \int\|x\|^s
K_\delta (x) \,dx \biggr]
\\
& = & 2^{s} \biggl[\int\|\theta\|^s \,dG(\theta) +
\delta^s \int\|x\|^s K(x) \,dx \biggr] < \infty,
\end{eqnarray*}
because $G$'s support points lie in a bounded subset
of $\real^d$.
Applying Lemma~\ref{Lem-ineq-moment}
to~(\ref{Eqn-bound-wasser-variational}), we obtain that for $\delta<
1$,
%
%e26 #&#
\begin{eqnarray}\label{Eqn-bound-L1}
\Wtwo^2\bigl(G*K_\delta,G'*K_\delta
\bigr) &\leq& C(d,K,s) \bigl\|G*K_\delta-G'*K_\delta
\bigr\|_{L_1}^{(s-2)/s}
\nonumber\\[-8pt]\\[-8pt]
& \leq& C(d,K,s) \bigl\|G*K_\delta-G'*K_\delta
\bigr\|_{L_2}^{2(s-2)/(d+2s)}.\nonumber
\end{eqnarray}
Here, constants $C(d,K,s)$ are different in each line, and they are dependent
only on $d,s$ and the $s$th moment of density function $K$.

Next, we use a known fact that for an arbitrary
(signed) measure $\mu$ on $\real^d$ and function
$g \in L_2(\real^d)$, there holds $\|\mu*g\|_{L_2} \leq|\mu| \|g\|_{L_2}$,
where $|\mu|$ denotes the total variation of $\mu$:
%
%e27 #&#
\begin{eqnarray}
\label{Eqn-bound-L2} \bigl\|G*K_\delta-G'*K_\delta
\bigr\|_{L_2} & = & \|p_G*g_\delta-
p_{G'}*g_\delta\|_{L_2}
\nonumber
\\
& = & \bigl\|(p_G-p_{G'})*g_\delta\bigr\|_{L_2}
\\
& \leq& 2\dV(p_G,p_{G'})\|g_\delta
\|_{L_2}.\nonumber
\end{eqnarray}
By Plancherel's identity,
\begin{eqnarray*}
\|g_\delta\|_{L_2}^2 & = & \frac{1}{(2\pi)^d}
\int\frac{\tildeK_\delta(\omega)^2} {
\tildef(\omega)^2} \,d\omega = \frac{1}{(2\pi)^d} \int_{\real^d}
\frac{\tildeK(\omega\delta
)^2}{\tildef(\omega)^2}\,d\omega
\\
& \leq& C \int_{[-1/\delta,1/\delta]^d} \tildef(\omega )^{-2}\,d\omega.
\end{eqnarray*}
The last bound holds because $\tildeK$ has support in $[-1,1]^d$ and
is bounded by a
constant.

Collecting equations (\ref{Eqn-perturb}), (\ref
{Eqn-bound-wasser-variational}), (\ref{Eqn-bound-L1}) and (\ref
{Eqn-bound-L2}) and the
preceding display, we have
\begin{eqnarray*}
&&\Wtwo^2\bigl(G,G'\bigr)
\\
&&\qquad\leq C(d,K,s) \biggl\{ \inf_{\delta\in(0,1)} \delta^2\\
&&\qquad\quad\hspace*{51pt}{} +
\dV(p_G,p_{G'})^{{2(s-2)}/({d+2s})} \\
&&\qquad\quad\hspace*{62pt}{}\times\biggl[\int
_{[-1/\delta,1/\delta]^d} \tildef(\omega)^{-2}\,d\omega
\biggr]^{({s-2})/({d+2s})} \biggr\}.
\end{eqnarray*}

If $|\tildef(\omega) \prod_{j=1}^{d}|\omega_j|^\beta| \geq d_0$ as
$\omega_j \rightarrow\infty$ $
(j=1,\ldots, d)$ for some positive constant~$d_0$, then
\begin{eqnarray*}
&&
\Wtwo^2\bigl(G,G'\bigr) \\
&&\qquad \leq C(d,K,s,\beta) \Bigl\{
\inf_{\delta\in
(0,1)} \delta^2 \\
&&\qquad\quad\hspace*{62.5pt}{} + \dV(p_G,p_{G'})^{{2(s-2)}/({d+2s})}
(1/\delta)^{{(2\beta+1)d(s-2)}/({d+2s})} \Bigr\}
\\
&&\qquad \leq C(d,K,s,\beta) \dV(p_G,p_{G'})^{
{4(s-2)}/({2(d+2s)+(2\beta+1)d(s-2)})}.
\end{eqnarray*}
The exponent tends to $4/(4+(2\beta+1)d)$ as $s\rightarrow\infty$,
so we obtain that
$\Wtwo^2(G,G') \leq C(d,\beta, r)\dV(p_G,p_{G'})^r$,
for any constant $r < 4/(4+(2\beta+1)d)$, as $\dV(p_G,p_{G'})
\rightarrow0$.

If $|\tildef(\omega) \prod_{j=1}^{d}\exp(|\omega_j|^{\beta})| \geq d_0$
as $\omega_j \rightarrow\infty$ $(j=1,\ldots, d)$
for some positive constants $\beta,d_0$, then
\begin{eqnarray*}
&&\Wtwo^2\bigl(G,G'\bigr)
\\
&&\qquad\leq C(d,K,s,\beta) \biggl\{ \inf_{\delta\in(0,1)} \delta^2 +
\dV(p_G,p_{G'})^{2(s-2)/(d+2s)} \exp-2d
\delta^{-\beta} \frac{s-2}{d+2s} \biggr\}.
\end{eqnarray*}
Taking $\delta^{-\beta} = -\frac{1}{d}\log\dV(p_G,p_{G'})$, we
obtain that
\[
d_{\rho^2}\bigl(G,G'\bigr) \leq C(d,\beta) \bigl(-\log
\dV(p_G,p_{G'})\bigr)^{-2/\beta}.
\]
\upqed\end{pf*}
\begin{pf*}{Proof of Lemma~\ref{Lem-gen-ineq}}
We exploit the variational characterization of $f$-divergences
(e.g.,~\cite{Nguyen-ieee-10}),
\[
\dPhi\bigl(f_i,f'_j\bigr) =
\sup_{\varphi_{ij}} \int\varphi_{ij} f'_j -
\phi^*(\varphi_{ij}) f_i \,d\mu,
\]
where the infimum is taken over all measurable functions defined on
$\Xcal$.
$\phi^*$ denotes the Legendre--Fenchel conjugate dual of convex
function $\phi$ [$\phi^*$ is again a convex function on $\real$
and is defined by $\phi^*(v) = \sup_{u\in\real}(uv - \phi(u))$].

%%Thus,
%%\[\drPhi(G,G') = \inf_{\vecq\in\Qcall} \sum_{ij}q_{ij}\sup_{
%%\int\varphi_{ij} f'_j - \phi^*(\varphi_{ij}) f_i.\]
%%
%For any $\qvec\in\Qcall$,
% = \sup_{\varphi}
%& = & \sup_{\varphi}
%= \sup_{\varphi}
%& \leq&
%= \sum_{ij} q_{ij} \sup_{\varphi_{ij}} \rho_\phi(f_i,f'_j),
%where the last inequality holds because the supremum is taken over a
%larger
%set of functions. Moreover, the bound holds for any $\vecq\in\Qcall$,
%so
%$\dPhi(p_G,p_{G'}) \leq\drPhi(G,G')$.
%}

By the variational characterization, $\rho_\phi$ is a convex
functional (jointly of its two arguments).
Thus, for any coupling $Q$ of two mixing measures $G$ and $G'$,
$\rho_\phi(p_G, p_{G'}) \hspace*{-0.3pt}=\hspace*{-0.3pt} \rho_\phi(\int\hspace*{-0.3pt} f(\cdot|\theta) \,dG,
\int\hspace*{-0.3pt} f(\cdot|\theta') \,dG') \hspace*{-0.3pt}=\hspace*{-0.3pt} \rho_\phi(\int f(\cdot|\theta) \,dQ,
\int\hspace*{-0.3pt} f(\cdot|\theta') \,dQ) \leq\break\int\rho_\phi(f(\cdot|\theta),
f(\cdot|\theta')) \,dQ$, where the inequality is obtained via
Jensen's inequality. Since this holds for any $Q$,
the desired bound follows.
\end{pf*}
\begin{pf*}{Proof of Proposition~\ref{Thm-1}}
(a) Suppose that the claim is not true, and
there is a sequence of $(G_{0},G) \in\Gcal_k(\Theta)\times\Gcal$
such that $\Wtwo(G_0,G_2) \geq r/2 > 0$ always holds and
that converges in $\Wtwo$ metric to $G_0^*\in\Gcal_k$ and
$G^*\in\Gcal$, respectively. This is due to the compactness
of both $\Gcal_k(\Theta)$ and $\Gcal$. We must have
$\Wtwo(G_0^*,G^*) \geq r/2 > 0$, so $G_0^* \neq G^*$.
%by Lemma~\ref{Lem-basic}.
At the same time,
$\dHel(p_{G_0^*},p_{G^*}) = 0$, which implies that $p_{G_0^*} =
p_{G^*}$ for almost all $x\in\Xcal$. By the finite identifiability condition,
$G_0^* = G^*$, which is a contradiction.

(b) is an immediate consequence of Theorem~\ref{Thm-bound-supnorm}, by
noting that under the given hypothesis, there is $c(k)>0$ depending
on $k$, such that
\begin{eqnarray*}
d_{h}^2(p_{G_0},p_{G}) & \geq&
\dV^2(p_{G_0},p_{G})/2
\\
& \geq& c(k,G_0) \Wtwo^4(G_0,G)
\end{eqnarray*}
for sufficiently small $\Wtwo(G_0,G)$. The boundedness
of $\Theta$ implies the boundedness of $\Wtwo(G_0,G)$, thereby extending
the claim for the entire admissible range of $\Wtwo(G_0,G)$.
(c) is obtained in a similar way to Theorem~\ref{Thm-convolution}.
\end{pf*}

%s5.2 #&#
\subsection{Proofs of posterior contraction theorems}

We outline in this section the proofs of Theorems
\ref{Thm-convergence-1} and~\ref{Thm-convergence-2}.
Our proof follows the same steps as in~\cite{Ghosal-Ghosh-vanderVaart-00},
with suitable modifications for the inclusion
of the Hellinger information function and the conditions
involving latent mixing measures.
The proof consists of results on the existence of tests, which are
turned into probability bounds on the posterior contraction.
\begin{pf*}{Proof of Lemma~\ref{Lem-convexball}}
We first consider case (I). Define $\Pcal_1 = \{p_{G}| G\in\Gcal\cap
B_W(G_1,r/2) \}$.
Since $\rho$ is a metric in $\Theta$,
it is a standard fact of Wasserstein metrics that $B_W(G_1,r/2)$
is a convex set.
Since $\Gcal$ is also convex, so is the set $\Gcal\cap B_W(G_1,r/2)$.
It follows that $\Pcal_1$ is a convex set of mixture distributions.
Now, applying a result from Birg\'e~\cite{Birge-84} and Le Cam (\cite
{LeCam-86}, Lemma 4, page 478),
there exist tests $\varphi_n$ that discriminate between
$P_{G_0}$ and convex set $\Pcal_1$ such that
%
%e28 #&#
%e29 #&#
\begin{eqnarray}
\label{Eqn-power-1} P_{G_0} \varphi_n & \leq& \exp
\bigl[-n \inf\dHel^2(P_{G_0},P_1)/2\bigr],
\\
\label{Eqn-power-2}
\sup_{G\in\Gcal\cap B_W(G_1,r/2)} P_{G}(1-\varphi_n) & \leq& \exp
\bigl[-n \inf\dHel^2(P_{G_0},P_1)/2\bigr],
\end{eqnarray}
where the exponent in the upper bounds are given by the infimum
Hellinger distance
among \textit{all} $P_1 \in\conv\Pcal_1 = \Pcal_1$.
Due to the triangle inequality, if
$\Wtwo(G_0, G_1) = r$ and $\Wtwo(G_1,G) \leq r/2$, then
$\Wtwo(G_0,G) \geq r/2$. So,
\[
\Psi_\Gcal(r) = \inf_{G\in\Gcal\dvtx  \Wtwo(G_0,G) \geq r/2} \dHel^2
(p_{G_0},p_{G}) \leq\inf\dHel^2(p_{G_0},P_1).
\]
This completes the proof of case (I).

Turning to case (II), for a constant $c_0 > 0$ to be determined,
consider a maximal $c_0r$-packing in the $\Wtwo$ metric in set $\Gcal
\cap B_W(G_1,r/2)$.
This yields a set of
$M(\Gcal,G_1,r) = D(c_o r, \Gcal\cap B_W(G_1,r/2), \Wtwo)$
points $\tilde{G}_1,\ldots, \tilde{G}_{M}$ in
$\Gcal\cap B_W(G_1, r/2)$. [In the following we drop the subscripts
of $M(\cdot)$.]

We note\vspace*{1pt} the following fact: For any $t=1,\ldots, M$,
if $G\in\Gcal\cap B_W(G_1,r/2)$ and
$\Wtwo(G,\tilde{G}_t) \leq c_0 r$, by Lemma~\ref{Lem-gen-ineq} we
have $\dHel^2(p_G,p_{\tilde{G}_t}) \leq d_{\rho_h^2}(G,\tilde{G}_t)
\leq\break C_1 d_{\rho^{2\alpha}}(G,\tilde{G}_t)
\leq  C_1 \Diam(\Theta)^{2(\alpha-1)}\Wtwo^2(G,\tilde{G}_t)
\leq  C_1 \Diam(\Theta)^{2(\alpha-1)} c_0^2r^2$
(the second inequality is due to the condition on the likelihood functions);
and so it follows that
\[
\dHel(p_{G_0},p_{G}) \geq\dHel(p_{G_0},p_{\tilde{G}_t})-
\dHel(p_{G},p_{\tilde{G}_t}) \geq\Psi_\Gcal(r)^{1/2}
- C_1^{1/2} \Diam(\Theta)^{\alpha-1} c_0r.
\]
Choose $c_0 = \frac{\Psi_\Gcal(r)^{1/2}}{ 2r\Diam(\Theta)^{\alpha
-1} C_1^{1/2}}$
so that the
previous bounds become
$\dHel(p_G,\break p_{\tilde{G}_t}) \leq\Psi_\Gcal(r)^{1/2}/2 \leq\dHel
(p_{G_0},p_{\tilde{G}_t})/2$ and
$\dHel(p_{G_0},p_{G}) \geq\Psi_\Gcal(r)^{1/2}/2$.

For each pair of $G_0,\tilde{G}_t$, there exist
tests $\omega_n^{(t)}$ of $p_{G_0}$ versus the closed Hellinger ball
$\{p_{G}\dvtx  \dHel(p_{G}, p_{\tilde{G}_t}) \leq
\dHel(p_{G_0},p_{\tilde{G}_t})/2 \}$ such that
\begin{eqnarray*}
P_{G_0} \omega_n^{(t)} & \leq& \exp\bigl[-n\dHel^2(P_{G_0},P_{\tilde{G}_t})/8\bigr],\\
\sup_{G\in\Gbar(\Theta): \dHel(p_G,p_{\tilde{G}_t}) \leq\dHel
(p_{G_0},p_{\tilde{G}_t})/2} P_{G}\bigl(1-\omega_n^{(t)}
\bigr) & \leq& \exp\bigl[-n \dHel^2(P_{G_0},P_{\tilde{G}_t})/8
\bigr].
\end{eqnarray*}
Consider the test $\varphi_n = \max_{t\leq M} \omega_n^{(t)}$, then
\begin{eqnarray*}
P_{G_0} \varphi_n & \leq& M \exp\bigl[-n
\Psi_\Gcal(r)/8\bigr],
\\
\sup_{G\in\Gcal\cap B_W(G_1,r/2)} P_{G}(1-\varphi_n) & \leq& \exp
\bigl[-n \Psi_\Gcal(r)/8\bigr].
\end{eqnarray*}
The first inequality is due to $\varphi_n \leq\sum_{t=1}^{M}\omega_n^{(t)}$, and
the second is due to the fact that for any $G\in\Gcal\cap B_W(G_1,r/2)$
there is some $t\leq M$ such that $\Wtwo(G,\tilde{G}_t) \leq c_0r$,
so that $\dHel(p_{G},p_{\tilde{G}_t}) \leq
\dHel(p_{G_0},p_{\tilde{G}_t})/2$.
\end{pf*}
\begin{pf*}{Proof of Lemma~\ref{Lem-ballcomplement}}
For a given $t \in\Nat$ choose a maximal $t\varepsilon/2$-packing for
set $S_t = \{G\dvtx  t\varepsilon< \Wtwo(G_0,G) \leq(t+1)\varepsilon\}$. This
yields a set $S'_t$ of at most $D(t\varepsilon/2,S_t,\Wtwo)$ points.
Moreover, every $G\in S_t$ is within distance $t\varepsilon/2$ of at least
one of the points in $S'_t$. For every such point
$G_1\in S'_t$, there exists a test
$\omega_n$ satisfying equations (\ref{Eqn-type-1}) and (\ref{Eqn-type-2}).
Take $\varphi_n$ to be the maximum of all tests attached this way to
some point $G_1\in S'_t$ for some $t \geq t_0$. Then, by the union bound
and the fact that $D(\varepsilon)$ is nonincreasing,
\begin{eqnarray*}
P_{G_0}\varphi_n &\leq&\sum_{t\geq t_0}
\sum_{G_1 \in S'_t} M(\Gcal,G_1,t\varepsilon) \exp
\bigl[-n \Psi_{\Gcal}(t\varepsilon)/8\bigr]
\\
&\leq& D(\varepsilon)\sum_{t\geq t_0}\exp\bigl[-n
\Psi_{\Gcal}(t\varepsilon)/8\bigr],
\\
\sup_{G\in\bigcup_{u\geq t_0}S_{u}} P_{G}(1-\varphi_n) &\leq&
\sup_{u\geq t_0} \exp\bigl[-n\Psi_{\Gcal}(u\varepsilon)/8\bigr] \leq \exp
\bigl[-n\Psi_{\Gcal}(t_0\varepsilon)/8\bigr],
\end{eqnarray*}
where the last inequality is due the monotonicity of $\Psi_{\Gcal
}(\cdot)$.
\end{pf*}
\begin{pf*}{Proof of Theorems~\ref{Thm-convergence-1} and
\ref{Thm-convergence-2}}
The proof for Theorem~\ref{Thm-convergence-1} proceeds in a
similar way to Theorem 2.1 of
\cite{Ghosal-Ghosh-vanderVaart-00}, while the proof for
Theorem~\ref{Thm-convergence-2} is similar to
their Theorem 2.4. The main difference is that the posterior
distribution statements are made with respect to mixing measure $G$
rather than mixture density $p_G$.
By a result of Ghosal et al.~\cite{Ghosal-Ghosh-vanderVaart-00} (Lemma
8.1, page 524),
for every $\varepsilon> 0$ and probability measure $\Pi$ on the set
$B_K(\varepsilon)$ defined by (\ref{Eqn-kl-neighborhood}), we have,
for every $C>0$,
\[
P_{G_0} \Biggl(\int\prod_{i=1}^{n}
\frac{p_G(X_i)}{p_{G_0}(X_i)} \,d\Pi(G) \leq\exp\bigl(-(1+C)n\varepsilon^2\bigr)
\Biggr) \leq\frac{1}{C^2n\varepsilon^2}.
\]
This entails that, for a fixed $C\geq1$, there is an event $A_n$
with $P_{G_0}$-probability at least $1-(C^2n\varepsilon_n^2)^{-1}$, for
which there holds
%
%e30 #&#
\begin{equation}
\label{Eqn-denominator} \int\prod_{i=1}^{n}p_G(X_i)/p_{G_0}(X_i)
\,d\Pi(G) \geq\exp\bigl(-2Cn\varepsilon_n^2\bigr)\Pi
\bigl(B_K(\varepsilon_n)\bigr).
\end{equation}
Let $\Ocal_n = \{G \in\Gbar(\Theta)\dvtx  \Wtwo(G_0,G) \geq M_n
\varepsilon_n\}$,
$S_{n,j} = \{G\in\Gcal_n\dvtx  \Wtwo(G_0,G)\in\break [j\varepsilon_n,
(j+1)\varepsilon_n)\}$
for each $j\geq1$.
The conditions specified by Lemma~\ref{Lem-ballcomplement} are satisfied
by setting $D(\varepsilon) = \exp(n\varepsilon_n^2)$ (constant in
$\varepsilon$).
Thus, there exist tests $\varphi_n$ for which equations (\ref{Eqn-type-1b})
and (\ref{Eqn-type-2b}) hold.
Then,
\begin{eqnarray*}
\hspace*{-3pt}&& P_{G_0}\Pi(G\in\Ocal_n|X_1,\ldots,X_n)
\\
\hspace*{-3pt}&&\qquad = P_{G_0}\bigl[\varphi_n \Pi(G\in
\Ocal_n|X_1,\ldots,X_n)\bigr] +
P_{G_0} \bigl[(1-\varphi_n)\Pi(G\in\Ocal_n|X_1,\ldots,X_n)\bigr]
\\
\hspace*{-3pt}&&\qquad \leq P_{G_0}\bigl[\varphi_n \Pi(G\in
\Ocal_n|X_1,\ldots,X_n)\bigr] +
P_{G_0}\indicator\bigl(A_n^c\bigr)
\\
\hspace*{-3pt}&&\qquad\quad{} + P_{G_0} \bigl[(1-\varphi_n)\Pi(G\in
\Ocal_n|X_1,\ldots,X_n)
\indicator(A_n)\bigr].
\end{eqnarray*}
%
%Exploiting Lemma~\ref{Lem-ballcomplement}, all terms in the
%preceeding display can be shown to vanish as $n\rightarrow\infty$.
%The proof for Theorem~\ref{Thm-convergence-1} proceeds in a
%similar way to Theorem 2.1 of
%Theorem~\ref{Thm-convergence-2} is similar to
%their Theorem 2.4.}

Due to Lemma~\ref{Lem-ballcomplement}, the first term in the preceding
display is bounded above by $P_{G_0}\varphi_n \leq
D(\varepsilon_n)\sum_{j\geq M_n} \exp[-n\Psi_{\Gcal_n}(j\varepsilon_n)/8]
\rightarrow0$, thanks to (\ref{Eqn-rate-cond}).
The second term in the above display is bounded by $(C^2n\varepsilon_n^2)^{-1}$ by
the definition of $A_n$. If
$n\varepsilon_n^2 \rightarrow\infty$, let $C=1$. If $n\varepsilon_n^2$
tends to a
positive constant away from 0, we let $C$ be arbitrarily large so that
this probability in the second term vanishes to 0.
To show that the third term in the display also vanishes as
$n\rightarrow\infty$,
we exploit the following expression:
\begin{eqnarray*}
&&
\Pi(G\in\Ocal_n|X_1,\ldots,X_n)
\\
&&\qquad= \int_{\Ocal_n}\prod_{i=1}^{n}p_G(X_i)/p_{G_0}(X_i)
\,d\Pi(G) \Big/ \int\prod_{i=1}^{n}p_G(X_i)/p_{G_0}(X_i)
\,d\Pi(G),
\end{eqnarray*}
and then obtain a lower bound for the denominator by (\ref
{Eqn-denominator}).
For the nominator, by Fubini's theorem,
%
%e31 #&#
\begin{eqnarray}\label{Eqn-T1}
&& P_{G_0} \int_{\Ocal_n \cap\Gcal_n} (1-\varphi_n)
\prod_{i=1}^{n} p_G(X_i)/p_{G_0}(X_i)
\,d\Pi(G)
\nonumber
\\
&&\qquad = P_{G_0} \sum_{j\geq M_n} \int
_{S_{n,j}}(1-\varphi_n)\prod
_{i=1}^{n} p_G(X_i)/p_{G_0}(X_i)
\,d\Pi(G)
\\
&&\qquad = \sum_{j\geq M_n} \int
_{S_{n,j}} P_{G}(1-\varphi_n) \,d\Pi(G)
\leq\sum_{j\geq M_n} \Pi(S_{n,j})\exp\bigl[-n
\Psi_{\Gcal_n}(j\varepsilon_n)/8\bigr],\nonumber
\end{eqnarray}
where the last inequality is due to (\ref{Eqn-type-2b}). In addition,
by (\ref{Eqn-support-1}),
%
%e32 #&#
\begin{eqnarray}\label{Eqn-T2}
&& P_{G_0} \int_{\Ocal_n \setminus\Gcal_n}(1-\varphi_n)
\prod_{i=1}^{n} p_G(X_i)/p_{G_0}(X_i)
\,d\Pi(G)
\nonumber
\\
&&\qquad = \int_{\Ocal_n \setminus\Gcal_n} P_G (1-\varphi_n)
\,d\Pi(G)
\\
&&\qquad \leq \Pi\bigl(\Gbar(\Theta) \setminus\Gcal_n
\bigr) = o\bigl(\exp\bigl(-2n\varepsilon_n^2\bigr)\Pi
\bigl(B_K(\varepsilon_n)\bigr)\bigr).\nonumber
\end{eqnarray}
Combining bounds (\ref{Eqn-T1}) and (\ref{Eqn-T2}) and condition
(\ref{Eqn-support-2}), we obtain
\begin{eqnarray*}
&& P_{G_0}(1-\varphi_n)\Pi(G\in\Ocal_n|X_1,\ldots,X_n)\indicator (A_n)
\\
&&\qquad\leq \frac{o(\exp(-2n\varepsilon_n^2)\Pi(B_K(\varepsilon_n))) + \sum_{j\geq M_n}
\Pi(S_{n,j})\exp[-n\Psi_{\Gcal_n}(j\varepsilon_n)/8]}{\exp
(-2n\varepsilon_n^2)\Pi(B_K(\varepsilon_n))}
\\
&&\qquad \leq o(1) + \exp\bigl(2n\varepsilon_n^2\bigr)\sum
_{j\geq M_n} \exp\bigl[-n\Psi_{\Gcal
_n}(j
\varepsilon_n)/16 \bigr].
\end{eqnarray*}
The upper bound in the preceding
display converges to 0 by (\ref{Eqn-rate-cond}), thereby
concluding the proof
of Theorem~\ref{Thm-convergence-2}. The proof of Theorem \ref
{Thm-convergence-1}
proceeds similarly.
\end{pf*}

%to the proof of Theorem~\ref{Thm-convergence-2} up
%to the application of \eqref{Eqn-type-2b}, where we now have:
%$P_{G_0} \int_{\Ocal_n \cap\Gcal_n} (1-\varphi_n)\prod_{i=1}^{n}
%p_G(X_i)/p_{G_0}(X_i) \Pi(G)
%and subsequently,
%{\exp(-2n\varepsilon_n^2)\Pi(B_K(\varepsilon_n))}\]
%which tends to 0 by conditions \eqref{Eqn-support-1s},
%and \eqref{Eqn-rate-conds-2s} of the theorem. The conclusion then
%follows.
%}

%s5.3 #&#
\subsection{Proof of other auxiliary lemmas}
\label{sec-aux}

\mbox{}

\begin{pf*}{Proof of Lemma~\ref{Lem-entropy}}
To simplify notation, we give a
proof for $W_1 \equiv d_\rho$. The general case for
$\Wr^r \equiv d_{\rho^r}$ can be carried out in the same way.

(a) Suppose that $(\eta_1,\ldots,\eta_T)$ forms an $\varepsilon$-covering
for $\Theta$ under metric $\rho$, where $T= N(\varepsilon,\Theta,
\rho)$ denotes the (minimum) covering number.
Take any discrete measure $G = \sum_{i=1}^{k}\p_i\delta_{\theta_i}$.
For each $\theta_i$ there is an approximating $\theta'_i$ among
the $\eta_j$'s such that
$\rho(\theta_i,\theta'_i) < \varepsilon$.
Let $\mathbf{p}'=(p'_1,\ldots,p'_k)$ be a $k$-dim vector in
the probability simplex
that deviates from $\pvec$ by less than $\delta$ in $l_1$ distance:
$\|\mathbf{p}'-\pvec\|_1 \leq\delta$. Define
$G'= \sum_{i=1}^{k} p'_i \delta_{\theta'_i}$. We shall argue that
\[
d_\rho\bigl(G,G'\bigr) \leq\sum
_{i=1}^{k} \bigl(p_i \wedge
p'_i\bigr)\rho\bigl(\theta_i,
\theta'_i\bigr) + \bigl\|\vecp- \vecp'
\bigr\|_1 \Diam(\Theta) \leq\varepsilon+ \delta\Diam (\Theta).
\]
[To see this, consider the following coupling between $G$ and $G'$:
associating $p_i \wedge p'_i$ probability mass of $\theta_i$ (from $G$)
with the same probability mass of $\theta'_i$ (from $G'$), while the
remaining mass from $G$ and $G'$ (of probability $\|\pvec-\mathbf
{p}'\|_1$)
are coupled in an arbitrary way. The right-hand side of the previous
display is an upper bound of the expectation of the $\rho$ distance
between two random variables distributed according to the described coupling.]

It follows that a $(\varepsilon+\delta\Diam(\Theta))$-covering for
$\Gcal_k(\Theta)$
can be constructed by combining each element of a $\delta$-covering
in the $l_1$ metric of the $k-1$-probability simplex and $k$ $\varepsilon
$-coverings
of $\Theta$. Now,\vspace*{1pt} the covering number of the $k-1$-probability simplex
is less
than the number of cubes of length $\delta/k$ covering $[0,1]^{k}$
times the volume of $\{(p'_1,\ldots,p'_k)\dvtx  p'_j \geq0,
\sum_{j}p'_j \leq1+\delta\}$, that is, $(k/\delta)^k(1+\delta)^{k}/k!
\sim(1+1/\delta)^k e^k/\sqrt{2\pi k}$.
It follows that $N(\varepsilon+\delta\Diam(\Theta),
\Gcal_k(\Theta),
d_\rho)
\leq T^k (1+1/\delta)^k e^k/\sqrt{2\pi k}$. Take $\delta= \varepsilon
/\Diam(\Theta)$
to\break achieve the claim.

(b) As before, let $(\eta_1,\ldots,\eta_T)$ be an $\varepsilon
$-covering of $\Theta$.
Take any $G = \sum_{i=1}^{k}\p_i\delta_{\theta_i}
\in\Gbar(\Theta)$, where $k$ may be infinity.
The collection of atoms $\theta_1,\ldots,\theta_k$
can be subdivided into disjoint subsets $S_1,\ldots, S_T$, some of which
may be empty, so that for each $t=1,\ldots, T$, $\rho(\theta_i,\eta_t)
\leq\varepsilon$ for any $\theta_i \in S_t$. Define
$\p'_t = \sum_{i=1}^{k}p_i \indicator(\theta_i\in S_t)$,
and let $G' = \sum_{t=1}^{T}\p'_t \delta_{\eta_t}$, then
we are guaranteed that
\[
d_\rho\bigl(G,G'\bigr) \leq\sum
_{i=1}^{k}\sum_{t=1}^{T}
p_i \indicator(\theta_i \in S_t) \rho(
\theta_i,\eta_{t}) \leq \varepsilon
\]
by using a similar coupling argument as in part (a).

Let $\mathbf{p}''=(p''_1,\ldots,p''_T)$ be a $T$-dim vector in
the probability simplex
that deviates from $\mathbf{p}'$ by less than $\delta$ in the $l_1$ distance:
$\|\mathbf{p}''-\mathbf{p}'\|_1 \leq\delta$. Take
$G''=\sum_{t=1}^{T}\p''_t \delta_{\eta_t}$. It is simple to observe
that $d_\rho(G',G'') \leq\Diam(\Theta)\delta$. By the triangle inequality,
\[
d_\rho\bigl(G,G''\bigr)\leq
d_\rho\bigl(G,G'\bigr) + d_\rho
\bigl(G',G''\bigr) \leq\varepsilon+ \delta
\Diam(\Theta).
\]

The foregoing arguments establish that an $(\varepsilon+\delta\Diam
(\Theta))$-covering
in the Wasserstein metric for $\Gbar(\Theta)$ can be constructed by combining
each element of the $\delta$-covering in $l_1$
of the $T-1$ simplex and a single covering of $\Theta$.
From the proof of part (a),
$N(\varepsilon+\delta\Diam(\Theta), \Gbar(\Theta), d_\rho) \leq
(1+1/\delta)^T e^T/\sqrt{2\pi T}$.
Take $\delta= \varepsilon/\Diam(\Theta)$ to conclude.

(c) Consider a $G = \sum_{i=1}^{k}p_i \delta_{\theta_i}$ such that
$d_\rho(G_0,G) \leq
2\varepsilon$. By definition, there is a coupling
$q \in\Qcal(\pvec,\mathbf{p}^*)$ so that $\sum_{ij} q_{ij}
\rho(\theta^*_i, \theta_j) \leq2\varepsilon$. Since $\sum_{j} q_{ij}
= p^*_i$,
this implies that $2\varepsilon\geq\sum_{i=1}^{k}p^*_i \min_{j} \rho
(\theta^*_i,\theta_j)$.
Thus, for each $i=1,\ldots,k$ there is a $j$ such that $\rho(\theta^*_i,\theta_j)
\leq2\varepsilon/p^*_i \leq2M\varepsilon$. Without loss of generality, assume
that $\rho(\theta^*_i,\theta_i) \leq2M\varepsilon$ for all
$i=1,\ldots,k$.
For sufficiently small $\varepsilon$, for any $i$, it is simple to
observe that
$d_\rho(G_0,G) \geq|p^*_i - p_i| \min_{j \neq i}\rho(\theta^*_i,
\theta_j)
\geq|p^*_i - p_i| \min_{j}\rho(\theta^*_i, \theta^*_j)/2$. Thus,
$|p^*_i - p_i| \leq4\varepsilon/m$.

Now, an $\varepsilon/4 + \delta\Diam(\Theta)$ covering in $d_\rho$ for
$\{G \in\Gcal_k(\Theta)\dvtx  d_\rho(G_0,G) \leq2\varepsilon\}$
can be constructed by combining the $\varepsilon/4$-covering for each
of the $k$ sets $\{\theta\in\Theta\dvtx  \rho(\theta,\theta^*_i) \leq
2M\varepsilon\}$
and the $\delta/k$-covering for each
of the $k$ sets $[p^*_i - 4\varepsilon/m,p^*_i+4\varepsilon/m]$. This entails
that
$N(\varepsilon/4+\delta\Diam(\Theta),
\{G \in\Gcal_k(\Theta)\dvtx  d_\rho(G_0,\break G) \leq2\varepsilon\}, d_\rho)
\leq[\sup_{\Theta'} N(\varepsilon/4,\Theta',\rho)]^k
(8\varepsilon k/m\delta)^k$.
Take $\delta= \varepsilon/(4\Diam(\Theta))$ to conclude the proof.
%
%(d) Any discrete measure $G \in\Gcal$ can be approximated by an
%element in
%$\Gcal_K(\Theta)$ with an approximation error $e^{-K\eta}$. That
%element
%is again approximated by another element in the $\varepsilon$-covering.
%Thus $N(e^{-K\eta}+\varepsilon, \Gcal, d_\rho) \leq N(\varepsilon, \Gcal_K(
%Pick $K = \log(1/\varepsilon)/\eta$ to achieve the desired bound.
%The proof of (e) is similar.
%}
\end{pf*}

%The above result shows that $\Gcal_k(\Theta)$'s are ``nice'' classes,
%for
%compactness and having a relatively small entropy number which grows
%linearly with respect to
%the number of atoms $k$ and the entropy of $\Theta$.
%$\Gcal(\Theta)$ is a noncompact set, while its closure set $\Gbar(
%entropy number that may grow exponentially with respect to the
%entropy of $\Theta$.
%This motivates us to consider more tractable
%approximating classes. A class of discrete measures $\Gcal\subset
%called $\varepsilon$-tight with respect to $\Gcal_k(\Theta)$ if
%with respect
%to $\Gcal_k(\Theta)$, then
%k(\log N(\varepsilon/2,\Theta,\rho)
%+ \log(e+2e\Diam(\Theta)/\varepsilon)).\]
%
%an element in
%$\Gcal_k(\Theta)$ with an approximation error $\varepsilon$. That element
%is again approximated by another element in the $\varepsilon$-covering.
%Thus $N(2\varepsilon, \Gcal, d_\rho) \leq N(\varepsilon, \Gcal_k(\Theta), d_
%Combining with Lemma~\ref{Lem-entropy}(a) to conclude.
%}
%
\begin{pf*}{Proof of Lemma~\ref{Lem-ineq-moment}}
(a) For arbitrary constant $R>0$,
we have $\int|p(x)-p'(x)|\|x\|^\kk \,dx \leq\int_{\|x\|\leq R}|p-p'|\|
x\|^\kk
+ \int_{\|x\|\geq R} (p+p')\|x\|^\kk
\leq R^\kk\|p-p'\|_{L_1} + \break R^{-(\s-\kk)}(\E_{p}\|X\|^\s+ \E_{p'}\|
X\|^\s)$, choosing
\[
R=\bigl[\bigl(\E_p\|X\|^\s+\E_{p'}\|X\|^\s
\bigr)/\bigl|p-p'\bigr|_{L_1}\bigr]^{1/\s}
\]
to conclude.

(b) For any $R>0$, we have
\begin{eqnarray*}
\int_{\|x\|\leq R} \bigl|p(x)-p'(x)\bigr|\,dx &\leq& V_d^{1/2}R^{d/2} \biggl[\int_{\|x\|
\leq R}
\bigl(p(x)-p'(x)\bigr)^2 \,dx\biggr]^{1/2} \\
&\leq& V_d^{1/2}R^{d/2}\bigl\|p-p'\bigr\|_{L_2}.
\end{eqnarray*}
We also have
\[
\int_{\|x\|\geq R} \bigl|p(x)-p'(x)\bigr|\,dx
\leq\int_{\|x\|\geq R} p(x)+p'(x) \,dx \leq R^{-s} \bigl(\E_p \|X\|^{s} +
\E_{p'}\|X\|^s\bigr).
\]
Thus,
\[
\bigl\|p-p'\bigr\|_{L_1} \leq\inf_{R>0}
V_d^{1/2}R^{d/2}\bigl\|p-p'
\bigr\|_{L_2} + R^{-s} \bigl(\E_p \|X
\|^{s} + \E_{p'}\|X\|^s\bigr),
\]
which gives the desired bound.
\end{pf*}

\section*{Acknowledgments}

The author wishes to thank Arash Amini, the Associate Editor and, in
particular, an anonymous reviewer for valuable comments and
suggestions. The shortened proof of Lemma~\ref{Lem-gen-ineq} was pointed out by the
reviewer, who also helped to spot a gap in a previous version of Lemma~\ref{Lem-convexball}.

%suskaldyti doi

% imsref loaded by lrinkeviciute, 2013-01-16 12:40:18
% imsref loaded by lrinkeviciute, 2013-01-16 12:45:21

\printaddresses

\end{document}